\documentstyle[12pt]{article}
\begin{document}

\newcommand{\de}{\mbox{$\delta $}}
\newcommand{\uA}{\mbox{$\underline{A}$}}
\newcommand{\uB}{\mbox{$\underline{B}$}}
\newcommand{\uS}{\mbox{$\underline{S}$}}
\newcommand{\e}{\mbox{${\bf \epsilon \rm}$}}
\newcommand{\qed}{\mbox{$ \qquad \qquad \qquad \qquad \qquad  \qquad \qquad \qquad \qquad
\qquad \qquad \qquad \qquad \qquad \qquad \qquad \qquad Q.E.D $ }}
\newtheorem{cor}{Corollary}[section]
\newtheorem{alg}{Algorithm}[section]
\newtheorem{lemma}{Lemma}[section]
\newtheorem{theo}{Theorem}[section]
\newtheorem{defi}{Definition}[section]
\newtheorem{exa}{Example} [section]
\newtheorem{pro}{Proposition}[section]
\newtheorem{rem}{Remark}[section]
\noindent {\bf THE MODEL MATCHING PROBLEM FOR A GENERAL CLASS OF NON-LINEAR DISCRETE INPUT - OUTPUT
SYSTEMS WITH CROSS-PRODUCTS. AN ALGORITHMIC APPROACH. \rm} \vskip 10 pt \noindent
{\it Stelios Kotsios, \par\noindent  University of Athens, Department of Economics,
\par\noindent
Division of Mathematics and Computer Science
\par\noindent
Pesmazoglou 8, 10559 Athens, Greece} \vskip 10 pt
\par\noindent
\small {\bf Abstract:} In this paper we examine 
the model matching problem that concerns
nonlinear input - output discrete systems, 
containing products among delays of input and output
signals,
through a special factorization. 
The
algebraic framework of $\de \e$-operators and the star-product that we adopt describe
these systems.
Moreover, a certain procedure that resembles the Euclidean division, allows
us to discover the linear factors of those systems, with respect to the above mentioned operations.
The entire approach is symbolically algorithmic and involves the use of suitable
software.

\vskip 5 pt \noindent {\bf Keywords:} Nonlinear, Discrete, Algebraic Geometry, Model Matching,  Algorithmic, Computational, Cross-products.
\normalsize
\section{Introduction}
As is well known, nonlinear systems are used in a variety  of applications and have
been the focus of research for a number of years. This, along with a recent increase
in their use, has led to a steady development of the relevant theory. Furthermore,
various advances in computer technology have had a major effect on control system
analysis and design. These advances have led to the development of new algorithms and
techniques, which in turn have allowed design methods to be accomplished with greater
speed and efficiency, \cite{kn:carcanias}.
 Yet, computational methods for the study of nonlinear systems are (still) at an early stage of development. Such computational methods have been applied in the study of modeling,
the problem of feedback linearization, and in the global optimization problem,
(\cite{kn:glad3},\cite{kn:fliessnew},\cite{kn:carcanias}, to mention but a few).
 All of these problems were investigated mostly in relation to continuous systems, and, in some cases, in relation to non-linear discrete
 time systems, \cite{kn:factn}.
By definition, a discrete time system evaluates input and output signals over a
countable number of time instants. In some cases, a discrete time system is obtained
from continuous time systems through sampling at certain time instants
\cite{kn:gladbook}. A further  case concerns  systems that are naturally and directly
described in discrete form, typically in financial or economic systems.
There is a rich literature devoted to the study of discrete systems. 
Certain works
approach the issue through analytical tools, as the Lie Derivatives and Distributions
\cite{kn:cyrot3},\cite{kn:monaco} 
and others by using algebraic methodologies, like
differential algebra or rings theory, 
\cite{kn:glad3},\cite{kn:fliessnew},\cite{kn:glad1},\cite{kn:perdon},\cite{kn:wu}.
Some researchers give algorithmic results \cite{kn:glad1} and others study the discrete systems
from the theoretical point of view \cite{kn:sontag},\cite{kn:sontag2}
An
important feature of discrete systems is their strong computational orientation. This
is connected to the fact that the procedures developed for these systems are more "
suitable " for realization in computational machines. Thus,  several design techniques
can be accomplished at a faster rate and in a way that makes them appropriate for
further applications. 
\par\noindent
The present paper focuses on nonlinear discrete input-output systems of the form:
\[  y(t)+ \sum a_i y(t-i)+\sum \sum a_{i_1i_2} y(t-i_1) y(t-i_2) + \cdots +\]
\[+\sum
\sum \cdots \sum a_{i_1i_2 \ldots i_n} y(t-i_1) \cdots y(t-i_n)=\]
\[=\sum b_j u(t-j) +\cdots +
\sum \sum \cdots \sum b_{j_1j_2 \ldots j_m} u(t-j_1) \cdots u(t-j_m)+ \]
\[+\sum\sum c_{kl}y(t-k)u(t-l)+\cdots +\sum \sum\sum \cdots \sum\sum
\cdots \sum \]
\begin{equation}\label{main}
                 c_{k_1k_2 \ldots k_{n'}l_1 l_2 \ldots
l_v} y(t-{k_1}) \cdots y(t-{k_{n'}}) u(t-l_1) \cdots u(t-l_v)
\end{equation}
Equation (\ref{main}) transforms causal input signals (i.e. $u(t)=0$ for $t<0$) to causal
output signals. A set of initial conditions
$y(0)=y_0$, $y(1)=y_1$, $\ldots$, $y(k)=y_k$ is always assigned to (\ref{main}) and the lowest
delay output term (that is $y(t)$) appears in the linear part of the system.
From all the above, we conclude that for any given input sequence the entire future output is uniquely determined,\cite{kn:kalouptsidis}.
The systems of the form (\ref{main}), which contain products among inputs and output signals, that are
sometimes called 'cross-products', encompass a broad variety of nonlinear discrete
systems. We obtain these either through transformations of nonlinear discrete
state-space representations into input - output forms \cite{kn:rugh}, or when we use
Taylor's expansion method to approximate other more general nonlinear discrete systems
\cite{kn:kalouptsidis},\cite{kn:rugh}. These are employed in signal processing theory,
whenever it is necessary to construct nonlinear representations of discrete signals,
(they are an extension of the infinite inpulse response filters to a nonlinear set up, 
\cite{kn:kalouptsidis}),
in nonlinear time-series analysis and in adaptive control, in the context of  designing
nonlinear adaptive controllers, \cite{kn:kokotovic}. Design problems for such systems
have been studied in the past through the use of either analytical or algebraic tools (see for
instance \cite{kn:glad3},\cite{kn:glad1},\cite{kn:kotta}). We have to mention here
especially the following works \cite{kn:sontag},\cite{kn:sontag2} which study the above systems 
by means of algebraic methods and provide us with theoretical results.
\par
\noindent
In this
paper we would like to propose a novel way of dealing with systems of the form
(\ref{main}) within a computer algebra environment. We do this initially by means of
special operators and operations that create an appropriate algebraic framework and
then through the development of symbolic computational algorithms that allow us to
solve the model matching's design problem.
Specifically, we use the so-called $\de \e$-operator to deal with cross-products. This
operator, introduced in \cite{kn:bibo}, is an extension of the simple $\de$ and
$\e$-operators that are used in the algebraic description of nonlinear discrete input
- output systems without cross-products \cite{kn:phd1}, and which are also an
extension of the simple shift operator $q$ in linear discrete systems
\cite{kn:amstrong}. Via these operators we can define the so-called $\de$,$\e$ and $\de \e$ - polynomilas.
These polynomials involve the delays of output and input signals in a nonlinear
polynomial way, and are
similar to polynomials with many variables. ( A special case of them are the linear polynomials ).
By means of these polynomials we can formally rewrite (\ref{main})
as follows: $ F[y(t),u(t)]=0$, where $F$ a $\de \e$-polynomial. Among these polynomials, we can define two
product operations: the dot-product, which corresponds to the usual product among
polynomials, and the star-product, which corresponds to the substitution of one
polynomial by another, or, in more system-oriented terminology, to the cascade
connection of systems. Though there is a similarity between this algebraic background
and others, there is also the interesting peculiarity that the set of $\de
\e$-polynomials with respect to the star-product is not a ring, \cite{kn:glad3}.
Indeed, the property
 $A \ast (B+C)=A\ast B+A \ast C$ is NOT valid for $\de \e$-polynomials. 
\par\noindent
The concept of formal $\de$ or $\e$-linear factors is presented here and departs from previous
approaches, \cite{kn:difeq}. The basic idea underlying this concept, is that it \"{} reveals \"{},
through appropriate algorithms, the linear $\de$ or $\e$-polynomials that are hidden
within a nonlinear $\de \e$-polynomial $F$ and they are factors of $F$ with respect either to the dot-product or to the star-product. Specifically, if $F$ is a given $\de
\e$-polynomial, then we can rewrite it as:
\begin{equation}\label{compactform}
F=\sum_{\bf \theta} \de_0^{\kappa_{\bf \theta}}\e_0^{\sigma_{\bf \theta}}\cdot (c_{\bf \theta}\de_{{\bf i}_{\bf \theta}}\e_{{\bf j}_{\theta}} \ast [L_{\bf \theta},M_{\bf \theta}])+{\bf R}_{\de,h}+{\bf R}_{\e,h}+{\bf R}_{\de\e}
\end{equation}
where, in general, the coefficients $c_{\bf \theta}$, the linear polynomials $L_{\bf \theta},M_{\bf \theta}$
and the remainders 
${\bf R}_{\de,h}$,${\bf R}_{\e,h}$,${\bf R}_{\de\e}$,
depend on the parameters $w_{ijhk},s_{ijhk}$.  These polynomials $L_{\bf \theta},M_{\bf \theta}$ are the linear formal factors of $F$, upon discussion. The immediate consequence of this factorization is that we
can transfer questions about the behaviour of the nonlinear system to questions about
the structure of the linear polynomials it \"{} contains \"{} and apply well-known
methods for dealing with those polynomials. The above expression of a given
$\de \e$-polynomial, is achieved by means of certain computational algorithms. These
algorithms are simply a kind of division among $\de \e$-polynomials with respect to
the star-product. Though this procedure resembles the method of Gr\"{o}bner basis
\cite{kn:cox}, there are some essential differences: (1) The set of $\de
\e$-polynomials is not a ring. (2) Here we are using  just one $\de \e$-polynomial,
instead of a set of such polynomials. (3) It involves the notion of free parameters.
We can then achieve certain tasks by giving
suitable values to those parameters. For instance, we can eliminate the remainders ${\bf R}_{\de,h},{\bf R}_{\e,h},{\bf R}_{\de\e}$ or we can find the common factors of the polynomials
$L_{\bf \theta},M_{\bf \theta}$, if any. All these can be implemented in a computer and an
appropriate software has already been created through the Mathematica
Package, which was used throughout this paper.
To elucidate the whole issue, let us present a short exhibitory example.
Let us have the bilinear system: $y(t)+y(t-1)=u(t-1)+y(t-2)u(t-1)$. This system can rewritten
as $(\de_0+\de_1-\e_1-\de_2 \e_1)[y(t),u(t)]=0$ or
$F[y(t),u(t)]=0$, where $F$ is the $\de$-polynomial $\de_0+\de_1-\e_1-\de_2 \e_1$.
$F$ can be written as:
\[ \de_0+\de_1-\e_1-\de_0\e_0\ast[w_{1,0,0,0}\de_0+w_{1,0,0,1}
\de_1+\de_2,s_{1,0,0,0}\e_0+\e_1]+s_{1,0,0,0}\de_0\e_0 \ast [w_{2,0,0,0}\de_0+\]
\[+w_{2,0,0,1}\de_1+\de_2,\e_0]+w_{1,0,0,1}\de_0\e_0 \ast [w_{3,0,0,0}\de_0+\de_1,s_{3,0,0,0}\e_0+\e_1]+(w_{1,0,0,1}s_{1,0,0,0}-\]
\[-w_{2,0,0,1}s_{1,0,0,0}-w_{1,0,0,1}s_{3,0,0,0})\de_0\e_0\ast [w_{4,0,0,0}\de_0+\de_1,\e_0]+\]
\[+(w_{1,0,0,0}-w_{1,0,0,1}w_{3,0,0,0})\de_0\e_0 \ast [\de_0,s_{5,0,0,0}\e_0+\e_1]+{\bf R}_{\de\e} \]
(we do not write ${\bf R}_{\de \e}$ explicitly, due to its large size). The exact meaning of the indexes in the above expression will be explained later. The set of the Formal Linear factors of $F$ is:
$\{ \de_0+\de_1,\e_1,$
$w_{1,0,0,0}\de_0+w_{1,0,0,1}
\de_1+\de_2,s_{1,0,0,0}\e_0+\e_1,$
$w_{2,0,0,0}\de_0+
w_{2,0,0,1}\de_1+\de_2,\e_0,$
$w_{3,0,0,0}\de_0+\de_1,s_{3,0,0,0}\e_0+\e_1$
$w_{4,0,0,0}\de_0+\de_1,$
$\de_0,s_{5,0,0,0}\e_0+\e_1\}$.
 By giving to the parameters $w_{i,0,0,j},s_{i,0,0,j}$
certain values we can take concrete expressions for the system.
For instance, if we put
$w_{1,0,0,0}=16$, $w_{1,0,0,1}=-10$, $w_{2,0,0,0}=1$, $w_{2,0,0,1}=1$, $w_{3,0,0,0}=-1$, $w_{4,0,0,0}=-\frac{33}{26}$,
$s_{1,0,0,0}=2$, $s_{3,0,0,0}=-3$, $s_{6,0,0,0}=-1$ we can eliminate the remainder ${\bf R}_{\de\e}$ and $F$ takes the 
form:
$F=\de_0+\de_1-\e_1$$-\de_0\e_0 \ast [16\de_0-10\de_1+\de_2,2\e_0+\e_1]+$
$2\de_0\e_0 \ast [\de_0+\de_1+\de_2,\e_0]-$
$10\de_0\e_0 \ast [-\de_0+\de_1, -3 \e_0+
\e_1] +$
$52\de_0\e_0 \ast [-\frac{33}{26}\de_0+
\de_1,\e_0]$$+6\de_0\e_0 \ast [\de_0, -\e_0+\e_1]$.
\par
\noindent
As a main application of the above methodology, we examine
the model-matching problem for those systems. This specific problem concerns the
discovery of linear feedback connections that force the original system to behave like
a desired linear one. An appropriate choice of such a desired system can turn the
problem into one that concerns  tracking. In order to clarify our way of facing up the problem, let us 
recall some well known ideas from the theory of linear systems. Let us have
the linear input - output system:
\[ a_0y(t)+a_1y(t-1)+ \cdots + a_ky(t-k)=b_1u(t-1)+b_2u(t-2)+ \cdots + b_mu(t-m) \]
we want to find a feedback connection so that the closed-loop system behaves like a given 
desired linear system. By means of the operator $q$, where $q^{-1}y(t)=y(t-1)$, we rewrite the original system
and the desired system algebraically as: $Ay(t)=Bu(t)$ , $A_dy^*(t)=B_du_c(t)$, where
$A,B,A_d,B_d$, proper polynomials of the single variable $q^{-1}$.
We define the feedback-law upon request to be of the form $Ru=Tu_c-Sy$, where $R,T,S$ 
polynomials of $q^{-1}$ to be determined.
Using the fact that in the linear case, the usual product among polynomials 
corresponds to the cascade connection of systems (i.e. if $y(t)=Gu(t)$ and $u(t)=Fv(t)$ then
$y(t)=GFv(t)$) we take the closed-loop system $(RA+BS)y(t)=BTu_c(t)$.
Assuming that $y=y^*$ we obtain the polynomials
$R,S,T$ by solving the equations $RA+BS=A_d$, $BT=B_d$, 
\cite{kn:amstrong}.
\par
\noindent
The question, which  naturally arises, is how can we extend the above approach in the 
case of nonlinear discrete input-output system. We firstly faced this problem in the paper
\cite{kn:phd1}, where we introduced the $\de$-operator in order to replace the $q$ operator 
and the star-product instead of the common product, (therefore the star-product corresponds 
to the cascade connection of nonlinear input-output systems). By means of those tools
we managed to solve the problem only for input-output systems with nonlinear output,
\underline{linear} input and \underline{without} cross-products between input and output signals.
Extension of those results can be found at \cite{kn:kotritt},\cite{kn:factn},\cite{kn:ecc05},\cite{kn:yamanaka}.
The aim of the current paper is to present a solution of the model matching problem for a general 
class of nonlinear input-output systems with nonlinearities in the output, in the input,
and \underline{with} cross-products. Essentially, what we are doing here is to reduce the
feedback designing questions to the linear factors $L_{\bf \theta},M_{\bf \theta}$, appeared at the original system. Indeed, if a feedback-law $u(t)=Sy(t)$ is available, then the system (\ref{main}) is transformed to a closed-loop form which is nothing else than a difference equation with respect to $y(t)$. the core of our approach is to find conditions which guarantee that all the terms of the said difference equation have a common linear factor. To achieve that we use $\de\e$-operators and some suitable algorithms to express (\ref{main}) as (\ref{compactform}). Then, we throw away the remainders $R_{\de\e}$,$R_{\e,h},R_{\de,h}$ and we seek for those values of the free parameters and the quantity $S$, which provide us with linear polynomials $L_\theta,M_\theta$, which are not prime among each other. Let us denote their common factor by $\Phi$. If the desired system coincides with $\Phi$ or it has common factoe with it, then the previous feedback-law solves the problem for proper initial conditions. Our empasis is on the fact that the original nonlinear design has been reduced to numerous  operations among linear polynomials, which can be carried out by menas of well known techniques of algebraic geometry.
\par
\noindent
Moreover, our approach  in the present
paper will be purely formal and we will not address any issues that pertain to
internal  stability.
We take the main advantages of our approach to be the following: (a) that it is a
symbolic algorithmic orientation that allows a direct computer implementation. (b)
that it allows the emergence of a set of solutions for each of the model matching problem above, as
opposed to a single solution. This is of vital importance, since it allows us to
choose between available solutions in a way that is informed by further constraints.
There is strong evidence that this approach, if successful, may be particularly useful
in the study of certain nonlinear problems within adaptive control or optimal
controllers design theory. A consideration of these problems, along with an
assessment  of the success of current methods on the basis of specific examples, will
be the subject matter of a future paper.
\par\noindent
The present paper is divided into two parts. In Part I we present the algebraic
framework and the algorithms, in Part II we are dealing with the model-matching
problem.
Throughout
the paper, ${\bf N,Z,Z^+,R}$ will denote the sets of natural, integers, positive integers, rational and real numbers,
respectively.

\section{The $\de$-operators and $\de$-polynomials}\label{prwto}
In this section we recall the notions of $\de$-operators and $\de$-polynomials. They have been appeared and studied in \cite{kn:phd1}, \cite{kn:kotritt}, \cite{kn:lappas1}. The scope of their presentation here is to constitute an introductory framework for the consideration of the $\de \e$-operators and $\de \e$-polynomials, which are going to be adopted at the next section. Let $y(t)$ be a real causal sequence, that is a function $y: {\bf Z} \to {\bf R}$
with the property $y(t)=0$, for $t<0$. We denote the set of such sequences by $F$. It is originating from the sampling of continuous functions. Let ${\bf i}$ be an ordered element of the set ${\bf Z^+}^n$, in other words, ${\bf i}=(i_1,i_2,\ldots,i_n)$ is a vector with $n$ positive integers as components, placed in an ascending way ($ i_1 \le i_2 \le \cdots \le i_n$). We call ${\bf i}$ a multi-index with dimension $n$. The $\de$-operator $\de_{\bf i}:F \to F$ is defined as follows:
\[\de_{\bf i}y(t)=\de_{i_1}\de_{i_2}\cdots\de_{i_n}y(t)=y(t-{i_1})y(t-{i_2})\cdots y(t-{i_n}) \]
If in the multi-index ${\bf i}$ certain components are equal among each other, then we can use an alternate notation. Indeed, let
$i_1=i_2=\cdots=i_a=j$, $i_{a+1}=i_{a+2}=\cdots=i_{a+b}=k$, $\ldots$,$i_{n-c}=i_{n-c+1}=\cdots=i_n=l$, then the operator $\de_{\bf i}$ can be rewritten as 
$\de_{\bf i}=\de_j^a\de_k^b\cdots \de_l^c$. The quantity $a+b+\cdots +c$ is called degree of $\de_{\bf i}$ and it is denoted by $\deg (\de_{\bf i})$. If ${\bf i}$ consists from one element ${\bf i}=(i)$, then we have the simple $i$-shift operator, i.e. $\de_iy(t)=y(t-i)$. If at least one of the components of a multi-index is equal to zero then the operator is called {\it zero} operator. A special case of zero operator is the identical operator $\de_0$, where $\de_0 y(t)=y(t)$. By convention we define the {\it null} operator, $\de_e$ as $\de_ey(t)=1$. In other words the symbol $e$ corresponds to the " empty " multi-index $( )$. The set of all $\de$-operators is denoted by $\Delta$.
\begin{exa} If ${\bf i}=(3)$, then $\de_{\bf i}y(t)=$
$\de_{3}y(t)=\de_3y(t)=y(t-3)$. If ${\bf i}=(1,1,1,2,3,3,4,4)$, then $\de_{\bf i}y(t)=$
$\de_{(1,1,1,2,3,3,4,4)}y(t)=\de_1^3\de_2\de_3^2\de_4^2y(t)=$
$y^3(t-1)y(t-2)y^2(t-3)y^2(t-4)$. The operators $\de_0y(t)=y(t)$, $\de_0^2\de_2^3y(t)=\de_{(0,0,2,2,2)}y(t)$$=y^2(t)y^3(t-2)$ are both zero operators.
\end{exa} 
Let ${\bf I}$ be set of multi-indexes. It may be ordered in a lexicographical way as follows: We say that the multi-index ${\bf i}=(i_1,i_2,\ldots,i_n)$ is " less " than the multi-index ${\bf j}=(j_1,j_2,\ldots,j_m)$, and we write ${\bf i} \prec {\bf j}$ if either $n<m$ or $n=m$ and the right-most nonzero entry of the vector ${\bf j} - {\bf i}$ is 
positive. For instance, $(10)\prec (9,9) \prec (1,1,3,3) \prec (1,2,3,3)$.
The order among multi-indexes implies an order among $\de$-operators in a natural way.
Indeed, we say that $\de_{\bf i} \prec \de_{\bf i'}$ if 
${\bf i} < {\bf i'}$. Therefore, $\de_{10} \prec \de_9^2 \prec \de_1^2\de_3^2
\prec \de_1\de_2\de_3^2$. We equip the
set of the $\de $-operators  $\Delta$, with two internal operations: the dot-product and the
star-product. The dot-product corresponds to the usual product among sequences, while
the star-product corresponds to the substitution of one sequence by another.
Specifically, let $y(t)$ be a causal sequence, ${\bf i}=(i_1,i_2,\ldots,i_n)$, ${\bf j}=(j_1,j_2,\ldots,j_m)$ two multi-indexes and $z(t)=\de_{\bf i}y(t)=$
$y(t-i_1)y(t-i_2)\cdots y(t-i_n)$, $w(t)=\de_{\bf j}y(t)=y(t-j_1)y(t-j_2) \cdots y(t-j_m)$. We define as dot-product of the operators $\de_{\bf i}$ and $\de_{\bf j}$ a new operator, denoted by 
$\de_{\bf i} \cdot \de_{\bf j}$, with the property: $\de_{\bf i}\cdot \de_{\bf j}y(t)=$$z(t)w(t)$$=y(t-i_1)y(t-i_2)\cdots y(t-i_n)$$y(t-j_1)y(t-j_2)\cdots y(t-j_m)$.
Their star-product is a new operator, denoted by 
$\de_{\bf i} \ast \de_{\bf j}$, with the property: $\de_{\bf i}\ast \de_{\bf j}y(t)=$$\de_{\bf i}w(t)$$=w(t-i_1)w(t-i_2)\cdots w(t-i_n)=$$y(t-j_1-i_1)y(t-j_2-i_1)\cdots y(t-j_m-i_1)\cdots$
$y(t-j_1-i_n)y(t-j_2-i_n)\cdots y(t-j_m-i_n)$. 
In order to give compact formulas for the star and dot-products, we need
the following operators among their multi-indexes.
Given two multi-indexes
$ {\bf i \rm}=(i_1,i_2,\ldots,i_k) $ and
$ {\bf j \rm}=(j_1,j_2,\ldots,j_{\lambda}) $, the new multi-index $ {\bf
i \oplus j \rm} $ is defined just juxtaposing
$ {\bf j \rm  }$ after $ {\bf i \rm} $.
Explicitly:
$ {\bf
i \oplus j \rm}=( i_1,j_1,i_2,i_3,j_2,\ldots,i_k,j_{\lambda}) $
where
$ i_1 \le j_1 \le i_2 \le i_3 \le j_2 \le \cdots \le i_k \le j_{\lambda} $.
We define the pointwise sum
$ {\bf j \rm} \dot + i $ as follows:
Let $ {\bf j \rm}=(j_1,j_2,\ldots,j_m) $ be a multi-index and $i$ an integer, then
$ {\bf j \rm}  \dot + i=(j_1+i,j_2+i,\ldots,j_m+i) $.
Using those notations we get the following formulaes:
\begin{pro}\label{formula}
The following properties are valid: 
(1) $ \de_{\bf i \rm} \cdot \de_{\bf j \rm}=\de_{\bf
i \oplus j \rm} $, (2) $ \de_{\bf i \rm} \ast \de_{\bf j \rm}=\de_{{\bf j}\dot +
i_1} \cdot \de_{{\bf j}\dot +
i_2} \cdots \de_{{\bf j}\dot +
i_n} $
\end{pro}
Further properties, concerning the dot and the star products, can be found
in \cite{kn:factn},\cite{kn:lappas1}.
We can also define an external operation, named addition, as follows:
$(\de_{\bf i}+\de_{\bf j})y(t)=\de_{\bf i}y(t)+\de_{\bf j}y(t)$.
It can be easily proved that $\de_{\bf k} \ast (\de_{\bf i}+\de_{\bf j})\ne$
$(\de_{\bf k} \ast \de_{\bf i})+(\de_{\bf k} \ast \de_{\bf j})$ and thus
the set
$ (  \Delta, +, \ast) $ is not a ring.
\par
\noindent
Expressions of the form $ A=\sum_{n=0}^w \sum_{\bf i \in I_n } a_{\bf i \rm} \de_{\bf i
\rm}$
are called $\de$-polynomials, where by ${\bf I_n \rm}$ we denote the set of multi-indexes with dimension $n$. By convention ${\bf
I_0 \rm}=\{\de_e\} $.
For each polynomial $ A$ we define
$ d(A)$ as follows:
$ d(A)=\min\{ \min ( i_1,i_2,\ldots,i_n)$,
${\bf i }=(i_1,i_2,\ldots,i_n) \in {\bf Z^+}^n$
 such that $a_{\bf i \rm} \ne 0 $, for $n=1,2,
\ldots,k \}$.
We define as degree of $A$ and we denote it by
$deg(A)$ the maximum $deg(\de)$ appeared in $A$.
The maximum term of a non-linear polynomial $A$, denoted by $max(A)$, is its largest term,
accordingly to the lexicographical order defined in the previous section.
 An expression of the form
$ \sum_{i \in {\bf Z}} a_i \de_i $ is called a linear polynomial.
\par \noindent
Two $\de$-polynomials $ A=\sum_{n=0}^\nu \sum_{\bf i \in I_n } a_{\bf i \rm} \de_{\bf i
}$ and $ B=\sum_{m=0}^\mu \sum_{\bf j \in J_m } b_{\bf j } \de_{\bf j
}$ are equal if $\nu=\mu$, ${\bf I}_n={\bf J}_m$, $n=0,1,\ldots, \nu$ and $a_{\bf i}=b_{\bf j}$.
Their dot-product is defined as follows:
$A\cdot B=\sum_{n=0}^\nu\sum_{m=0}^\mu\sum_{\bf i \in I_n }\sum_{\bf j \in J_m }a_{\bf i \rm}
 b_{\bf j }\de_{\bf i
}\cdot \de_{\bf j
}$.
This dot-product is nothing else than the classical product among polynomials with many
variables. In order to define their star-product, we proceed as follows: Let us consider the $\de$-polynomials $A$ and $B$ as functions, transforming sequences to sequences, accordingly to the rules:
$A: F \to F, w(t) \to Aw(t) $,
$B: F \to F, y(t) \to By(t) $,
then, the star-product $A \ast B$ is defined as:
\[ A \ast B: F \to F \quad,\quad y(t) \to A \ast B y(t)=A \circ B y(t)=A(By(t)) \]
There are certain formulaes for the calculation of the star-product. A first one appeared in \cite{kn:phd1}, another one is provided in \cite{kn:lappas1}. For the sake of completeness of our presentation we shall present the first one, which is the more frequently used. 
\begin{pro}
Let us suppose that we have the $\de$-polynomials $A=\sum_{n=0}^\nu\sum_{{\bf i}=(i_1,i_2,\ldots,i_n) \in {\bf I}_n} a_{\bf i}\de_{\bf i} $,
$B=\sum_{m=0}^\mu\sum_{{\bf j} \in {\bf J}_m} b_{\bf j}\de_{\bf j} $, their star-product is given by the relation:
\[ A \ast B=\sum_{n=0}^\nu\sum_{{\bf i}=(i_1,i_2,\ldots,i_n) \in {\bf I}_n}\sum_{({\bf j}_1,{\bf j}_2,\ldots,{\bf j}_n) \in (\cup_m {\bf J}_m)^n}a_{\bf i}b_{{\bf j}_1}b_{{\bf j}_2}\cdots b_{{\bf j}_n}\de_{{\bf j}_1 \dot{+} i_1}\de_{{\bf j}_2 \dot{+} i_2}\cdots \de_{{\bf j}_n \dot{+} i_n}\]

\end{pro}
Many times, for the calculation of simple star-products, we use the following properties:
\begin{pro}\label{propo0}
(1) $\de_{i}\ast A=\sum_{n=0}^\nu \sum_{\bf i \in I_n } a_{\bf i \rm} \de_{{\bf i
}\dot{+}i}$, $i \in {\bf Z}$ $\quad$ 
(2) $ \de_{\bf i} \ast A= ( \de_{i_1} \ast A ) \cdot  ( \de_{i_2} \ast A )  \cdots   (
\de_{i_n} \ast A ) $, ${\bf i}={i_1,i_2,\ldots,i_n})$ $\quad$
(3) $ \de_{\bf i} \ast (A \cdot B)= (\de_{\bf i} \ast A) \cdot (\de_{\bf i} \ast B) $
\end{pro}
\noindent
The following propositions, can be found in \cite{kn:factn}.
\begin{pro}\label{propo1}
$ (1)  [A+B]\ast C=A \ast C+B \ast C $,
$ (2) C \ast [A+B]=C \ast A+C \ast B$ iff $C$ linear,
$ (3) A \ast B \ne B \ast A $,
$ (4)  (A \ast B) \ast \Gamma= A \ast ( B \ast \Gamma) $,
$ (5) d(A \ast B)=d(A)+d(B) $,
$ (6) deg(A \ast B)=deg(A) \cdot deg(B) $
$(7)  (A \cdot B) \ast C=(A \ast C) \cdot (B \ast C) $.
\end{pro}
The following property is very useful. It claims that working with linear $\de$-polynomials and the star-product it is like working with polynomials of a single variable and the classical product among them.  
\begin{pro}
Let ${\bf L}$ be the set of linear $\de$-polynomials, then the set
$({\bf L},\ast,+)$ is a commutative ring and it is isomorphic to the ring $(R[x],\cdot,+)$, where $R[x]$ is the set of real polynomials o a single variable and $\cdot$ the operation of the polynomial product.
\end{pro}
{\bf Proof:} That $({\bf L},\ast,+)$ is a ring comes as a straightforward result of the proposition \ref{propo1}. For any linear $\de$-polynomial $M=\sum_{i=0}^km_i\de_i$ we define the map
\[ \varphi : ({\bf L},\ast,+)\to (R[x],\cdot,+) \quad , \quad \varphi(M)=\sum_{i=0}^km_ix^i\]
It can be easily proved that $\varphi(M+N)=\varphi(M)+\varphi(N)$, $N$ another linear $\de$-polynomail and that $\varphi$ is one-to-one and onto. Let now $N=\sum_{j=0}^h n_j\de_j$, then
$M \ast N=\sum_{i=0}^k\sum_{j=0}^h m_i n_j\de_{i+j}$. This means that $\varphi(M\ast N)=\sum_{i=0}^k\sum_{j=0}^h m_i n_j x^{i+j}=\sum_{\theta=0}^{k+h}(\sum_{s=0}^\theta m_s n_{\theta-s})x^{\theta} $
$=(\sum_{i=0}^k m_i x^i)\cdot(\sum_{j=0}^h n_j x^j)=\varphi(M) \cdot \varphi(N)$. This relation ensures that $\phi$ is an isomorphism and the theorem has been proved.
\section{The $\de \e$-operators and the $\de\e$-polynomials}
The $\de$-operators defined above, act on single sequences. This concept can help us to describe the delays appeared either in the input or in the output signal only. Nevertheless, a problem arises when we have to deal with cross-products. Specifically, let us suppose for instance, that we have the product $y^2(t-1)u^3(t-2)$, then using $\de$-operators we get $\de_1^2y(t)\de_2^3u(t)$. Hereafter, it is not clear how we shall handle this expression. The main difficulty is that we do not know which operator acts where. To pass over this obstacle we have to modify the notions of $\de$-operators and $\de$-polynomials in a proper way. Actually, we introduce  the $\de\e$-operator which acts instead of to a single sequence to a pair of sequences. This initially appeared in
\cite{kn:bibo} and has been examined in \cite{kn:lappas1}.
Let $y(t),u(t)$ be two causal, real sequences, defined over the set of positive integers. Let $i,j$ be integers. We define the $\de_i \times \de_j$-operator as an operator acting on the pair $[y(t),u(t)]$ (actually, we had to denote this pair by $(y(t),u(t))$, but we use brackets in order to avoid so many parenthesis), as follows:
\[ \de_i \times \de_j:  F\times F \to F : \de_i \times \de_j [y(t),u(t)]=y(t-i)u(t-j) \] 
Let ${\bf i}=(i_1, i_2, \ldots, i_m), {\bf j }=(j_1, j_2, \ldots, j_n) $ be multi-indexes.
The operator 
$ \de_{\bf i} \times \de_{\bf j} : F \times F \mapsto F $ is defined as:
\[ \de_{\bf i} \times \de_{\bf j}{[y(t), u(t)]}= y(t-i_1)y(t-i_2) \cdots
y(t-i_m)u(t-j_1)u(t-j_2) \cdots u(t-j_n) \]
This means that the operator $\de_{\bf i}$
acts exclusively on " outputs " and $\de_{\bf j}$ exclusively on " inputs ".
Sometimes, for the sake of appearance, the following  notation may be more convenient: $ \de_{\bf i} \times
\de_{\bf j}=\de_{\bf i} \epsilon_{\bf j} $. Therefore, the $\e$-operator is just an operator with poperties identical similar to the properties of $\de$-operator, except that it acts only on the second sequence (input). We call these operators $ \de \e
$-operators and we denote their set by ${\cal D }$. 
Obviously, $ \de_{\bf i} \times \de_e
[y(t), u(t)]=\de_{\bf i}y(t) $, $ \de_e \times \de_{\bf j} (y, u)=\e_{\bf j}u(t) $.
 A $\de\e$-operator, $\de_{\bf i}\e_{\bf j}$
with the property that $i_1=j_1=0$, (in other words, the lowest delays of the $\de$ and $\e$-parts of $\de_{\bf i}\e_{\bf j}$ are zero), is called a {\it zero} $\de\e$-operator. A special case of a zero $\de\e$-operator is the operator $\de_0\e_0$, with $\de_0\e_0[y(t),u(t)]=y(t)u(t)$. 
\begin{exa}
Let ${\bf i}=(0,1,1,2)$, ${\bf j}=(1,1,3,3,3)$, ${\bf r}=(0,0,1,2)$ and ${\bf
s}=(0,1,1,2)$, then $ \de_{\bf i} \e_{\bf j}[y(t),u(t)]=$ $\de_{(0,1,1,2)}$
$\e_{(1,1,3,3,3)}[y(t),u(t)]=$ $\de_0\de_1\de_1\de_2$ $\e_1\e_1\e_3\e_3\e_3
[y(t),u(t)]$ $=\de_0\de_1^2\de_2$ $\e_1^2\e_3^3[y(t),u(t)]=$
$y(t)y^2(t-1)y(t-2)u^2(t-1)u^3(t-3)$, furthermore $\de_{\bf
r}\e_{e}[y(t),u(t)]=\de_{\bf r}y(t)=$ $\de_0^2\de_1\de_2y(t)=y^2(t)y(t-1)y(t-2)$ and
$\de_e \e_{\bf s}[y(t),u(t)]=\e_{\bf s}u(t)=$ $\e_0\e_1^2\e_2 u(t)=u(t)u^2(t-1)u(t-2)$
\end{exa}
The order among multi-indexes implies an order among $\de \e$-operators in a natural way.
Indeed, we say that $\de_{\bf i}\e_{\bf j} \preceq \de_{\bf i'}\e_{\bf j'}$ if either
${\bf i} < {\bf i'}$ or ${\bf i} = {\bf i'}$ and ${\bf j} < {\bf j'}$. We equip the
set of the $\de \e$-operators with two internal operations: the dot-product and the
star-product.
Let $\de_{\bf i} \times \de_{\bf j}=\de_{\bf i}\e_{\bf j}$, 
$\de_{\bf i'} \times \de_{\bf j'}=\de_{\bf i'}\e_{\bf j'}$
be two $\de\e$-operators. As their dot-product, denoted
by $\de_{\bf i}\e_{\bf j}\cdot \de_{\bf i'}\e_{\bf j'}$, we define the operator:
\[ ( \de_{\bf i} \cdot \de_{\bf i'}) \times (\de_{\bf j} \cdot \de_{\bf j'})=
( \de_{\bf i} \cdot \de_{\bf i'})  (\e_{\bf j} \cdot \e_{\bf j'})=\de_{\bf i
\oplus i'} \e_{\bf j \oplus j'}\]
It can be easily proved that the dot-product corresponds to the usual product among sequences.
\begin{exa} Let ${\bf i}=(0,1)$, ${\bf j}=(1,1,2)$, ${\bf i'}=(0,0,1)$ and ${\bf j'}=(0,0)$.
Then, $\de_{\bf i} \e_{\bf j}$ $ \cdot \de_{\bf i'} \e_{\bf j'}[y(t),u(t)]=$
$\de_0\de_1 \e_1^2\e_2$ $ \de_0^2\de_1\e_0^2 [y(t),u(t)]=$
$\de_0^3\de_1^2\e_0^2\e_1^2\e_2 [y(t),u(t)]$ $=y^3(t)y(t-1)u^2(t)u^3(t-1)u(t-2)$.
Since ${\bf i} \oplus {\bf i'}=(0,0,0,1,1)$ and $ {\bf j} \oplus {\bf j'}=(0,0,1,1,2)$
we can easily verify the relevant formula.
\end{exa}
\noindent
As we mentioned earlier, when we were dealing with $\de$-operators, the star-product corresponded to the substitution of one sequence into another. In the case of $\de\e$-operators the star-product corresponds to the substitution of two sequences with cross-products into the "$y$-part" and the "$u$-part" respectively, of a given sequence. To clarify this concept, let us see an example. Let
$z(t)=y(t-1)y(t-2)u^2(t-1)$$=\de_1\de_2\e_1^2[y(t),u(t)]$ be a sequence, involving cross-products of the delays of the sequences $y(t)$ and $u(t)$. Moreover, we suppose that each of them is a cross-product of the delays of two other sequences $w(t),v(t)$:
$y(t)=w^2(t-2)v(t-3)=\de_2^2\e_3[w(t),v(t)]$,
$u(t)=w(t-1)v(t-1)v(t-2)=\de_1\e_1\e_2[w(t),v(t)]$.
By executing the substitution we get:
$z(t)=w^2(t-2)w^2(t-3)w^2(t-4)$$v(t-2)v(t-3)v(t-4)v(t-5)$
$=\de_2^2\de_3^2\de_4^2\e_2\e_3\e_4\e_5[w(t),v(t)]$, this latter $\de\e$-operator is the star-product of the operator $\de_1\de_2\e_1^2$ and the pair $(\de_2^2\e_3, \de_1\e_1\e_2)$. 
Let us present the whole subject formally. We have the $\de\e$-operators
$\de_{\bf i}\e_{\bf j}$,
$\de_{{\bf i}_1}\e_{{\bf j}_1}$,
$\de_{{\bf i}_2}\e_{{\bf j}_2}$. Let $[w(t),v(t)] \in F \times F$ be an arbitrary given pair
of causal sequences. We define the maps:
$\Phi_1: [w(t),v(t)]\to y(t)= \de_{{\bf i}_1}\e_{{\bf j}_1}[w(t),v(t)] $, $ \Phi_2: [w(t),v(t)]\to u(t)= \de_{{\bf i}_2}\e_{{\bf j}_2}[w(t),v(t)] $,
$ \Phi : [y(t),u(t)]\to z(t)= \de_{{\bf i}}\e_{{\bf j}}[y(t),u(t)] $.
The star-product of the operator $\de_{\bf i}\e_{\bf j}$ and the pair of operators
$[\de_{{\bf i}_1}\e_{{\bf j}_1},\de_{{\bf i}_2}\e_{{\bf j}_2}]$ is the $\de\e$-operator $\de_{\bf h}\e_{\bf s}$, which transfers the pair $[w(t),v(t)]$ to the sequence $z(t)$. In other words, it corresponds to the map:
\[\Phi \circ ( \Phi_1, \Phi_2): [w(t),v(t)]\to z(t)=\de_{\bf h}\e_{\bf s}[w(t),v(t)] \]
Alternatively, the star-product is denoted by
$\de_{\bf i}\e_{\bf j} \ast [\de_{{\bf i}_1}\e_{{\bf j}_1},\de_{{\bf i}_2}\e_{{\bf j}_2}]$.
Of course, by using the null operators $\de_e$, $\e_e$ we can include the star-product among simple $\de$-operators and $\de \e$-operators into the above definition. For instance:
$\de_{\bf i}\e_e \ast [\de_{{\bf i}_1}\e_{{\bf j}_1},
\de_{e}\e_{e}]=\de_{\bf i} \ast \de_{{\bf i}_1}\e_{{\bf j}_1}$.
Let us equip the set of $\de\e$-operators with an external operation, named addition, as follows:
\[ [\de_{{\bf i}_1}\e_{{\bf j}_1}+
\de_{{\bf i}_2}\e_{{\bf j}_2}][y(t),u(t)]=
\de_{{\bf i}_1}\e_{{\bf j}_1}[y(t),u(t)]+
\de_{{\bf i}_2}\e_{{\bf j}_2}][y(t),u(t)]\]
The following properties are valid for all the above operations. Their
proofs are either obvious from the definitions and therefore omitted, or they can be found in \cite{kn:bibo},\cite{kn:lappas1}.
\begin{pro}\label{proti}
(a) The star-product is well defined.
\par
\noindent
(b) $ \de_{\bf i}\e_{\bf j} \ast [\de_{{\bf i}_1}\e_{{\bf j}_1},
\de_{{\bf i}_2}\e_{{\bf j}_2}]=$ $\de_{{\bf i}_1\dot{+}i_1} \de_{{\bf i}_1\dot{+}i_2} \cdots \de_{{\bf
i}_1 \dot{+}i_n} $ 
$\e_{{\bf j}_1\dot{+}i_1} \e_{{\bf j}_1\dot{+}i_2} \cdots \e_{{\bf
j}_1 \dot{+}i_n} $ 
$\de_{{\bf i}_2\dot{+}j_1} \de_{{\bf i}_2\dot{+}j_2} \cdots \newline \de_{{\bf
i}_2 \dot{+}j_n} $ 
$\e_{{\bf j}_2\dot{+}j_1} \e_{{\bf j}_2\dot{+}j_2} \cdots \e_{{\bf
j}_2 \dot{+}j_n} $
where ${\bf i}=(i_1,i_2,\ldots,i_n)$ and ${\bf j}=(j_1,j_2,\ldots,j_n)$ 
\par
\noindent
(c) $\de_{\bf i} \ast (\de_{{\bf i}_1}\e_{{\bf j}_1})=$
$\de_{\bf i} \ast (\de_{{\bf i}_1} \times \de_{{\bf j}_1})=$
$(\de_{\bf i} \ast \de_{{\bf i}_1}) \times (\de_{\bf i} \ast \de_{{\bf j}_1})=$
$(\de_{\bf i} \ast \de_{{\bf i}_1})(\e_{\bf i} \ast \e_{{\bf j}_1})$
\vskip 10 pt
\noindent
(d) $ \de_{\bf i}\e_{\bf j} \ast [\de_{{\bf i}_1}\e_{{\bf j}_1}+
\de_{{\bf i}_2}\e_{{\bf j}_2}]\ne$
$ \de_{\bf i}\e_{\bf j} \ast \de_{{\bf i}_1}\e_{{\bf j}_1}+
\de_{\bf i}\e_{\bf j} \ast\de_{{\bf i}_2}\e_{{\bf j}_2}$
\end{pro}
\begin{exa}
Let ${\bf i}=(0,1)$, ${\bf j}=(1,1,2)$,
${\bf i}_1=(0,0,1)$, ${\bf j}_1=(1,1)$,
${\bf i}_2=(0,1)$, ${\bf j}_2=(2)$, then
$\de_{\bf i}\e_{\bf j} \ast [\de_{{\bf i}_1}\e_{{\bf j}_1},
\de_{{\bf i}_2}\e_{{\bf j}_2}]$
$=\de_0\de_1\e_1^2\e_2 \ast [\de_0^2\de_1\e_1^2,\de_0\de_1\e_2]$. By setting
$y(t)=\de_0^2\de_1\e_1^2[w(t),v(t)]=w^2(t)w(t-1)v^2(t-1)$,
$u(t)=\de_0\de_1\e_2[w(t),v(t)]=w(t)w(t-1)v(t-2)$ and
$z(t)=\de_0\de_1\e_1^2\e_2[y(t),u(t)]=$
$y(t)y(t-1)u^2(t-1)u(t-2)$,
we finally get
$z(t)=w^2(t)w^5(t-1)w^4(t-2)w(t-3)$
$v^2(t-1)v^2(t-2)v^2(t-3)v(t-4)=$
$\de_0^2\de_1^5\de_2^4\de_3$
$\e_1^2\e_2^2\e_3^2\e_4[w(t),v(t)]$
and thus
$ \de_0\de_1\e_1^2\e_2 \ast [\de_0^2\de_1\e_1^2,\de_0\de_1\e_2]=\de_0^2\de_1^5\de_2^4\de_3
\e_1^2\e_2^2\e_3^2\e_4$. 
Using  part (b) of proposition (\ref{proti}) we get:
$ \de_{(0,1)}\e_{(1,1,2)}\ast$ $ [ \de_{(0,0,1)}\e_{(1,1)},\de_{(0,1)}\e_2]$ $=\de_{(0+0,0+0,1+0)}\de_{(0+1,0+1,1+1)}$ $\e_{(1+0,1+0)}
\e_{(1+1,1+1)}$$\de_{(0+1,1+1)}\de_{0+1,1+1)}$$\de_{(0+2,1+2)}\newline\e_{2+1}\e_{2+1}\e_{2+2}$
$=\de_0^2\de_1^5\de_2^4\de_3
\e_1^2\e_2^2\e_3^2\e_4$
\end{exa}
\noindent
The $\de \e$-polynomials are straightforward extensions of the $\de \e$-operators that
are necessary for the description of discrete non-linear systems with products among
input and output sequences. These have been introduced and studied in
\cite{kn:bibo},\cite{kn:lappas1}. We will go over the main points, adding some new
results that will clarify the ideas in the present paper.
Let ${\bf I}_n$, ${\bf J}_m$ be sets of multi-indexes with dimensions $n$ and $m$ respectively. By convention we set ${\bf I}_0={\bf J}_0=\{e\}$. An expression of the form:
\[ A=\sum_{n=0}^\nu\sum_{m=0}^\mu\sum_{({\bf i,j}) \in {\bf I}_n \times {\bf J}_m}c_{\bf ij}\de_{\bf i}\e_{\bf j} \]
where $c_{\bf ij}$ are real numbers, is called a $\de\e$-polynomial. Whenever a $\de\e$-polynomial acts on a pair of sequences, produces a polynomial of delays of those sequences, including also and products between delays of the first and second sequence, (the so-called cross-products). By means of the null index $e$, we can decompose a $\de\e$-polynomial into its pure $\de$-part, $\e$-part and $\de\e$-part. Indeed,
\[ A_{\de}=\sum_{n=0}^\nu\sum_{({\bf i,j}) \in {\bf I}_n \times {\bf J}_0}c_{\bf ij}\de_{\bf i}\e_{\bf j}= \sum_{n=0}^\nu\sum_{{\bf i} \in {\bf I}_n }c_{\bf ie}\de_{\bf i}\e_{e}=
\sum_{n=0}^\nu\sum_{{\bf i} \in {\bf I}_n }c_{\bf i}\de_{\bf i}\]
is the pure $\de$-part of $A$ and correspondingly
\[ A_{\e}=\sum_{m=0}^\mu\sum_{({\bf i,j}) \in {\bf I}_0 \times {\bf J}_m}c_{\bf ij}\de_{\bf i}\e_{\bf j}= \sum_{m=0}^\mu\sum_{{\bf j} \in {\bf J}_m }c_{\bf ej}\de_{e}\e_{\bf j}=
\sum_{m=0}^\mu\sum_{{\bf j} \in {\bf J}_m }c_{\bf j}\e_{\bf j}\]
the pure $\e$-part. Actually, it is a $\de$-polynomial acting only to inputs.
The quantity $A_{\de\e}=A-A_{\de}-A_{\e}$ is the pure $\de\e$-part. We shall call expressions of the form:
\begin{equation}\label{l}
\sum_{(i,j) \in ({\bf I}_1 \times {\bf J}_0) \cup ({\bf I}_0 \times {\bf J}_1)}c_{ij}\de_i\e_j=
\sum_{i \in {\bf I}_1}c_i\de_i+\sum_{j \in {\bf J}_1} c_j \e_j 
\end{equation}
\[ \sum_{(i,j) \in ({\bf I}_1 \times {\bf J}_0)}c_{ie}\de_i\e_e=\sum_{i \in {\bf I}_1}c_i\de_i \quad,\quad
\sum_{(i,j) \in ({\bf I}_0 \times {\bf J}_1)}c_{ej}\de_e\e_j=\sum_{j \in {\bf J}_1}c_j\e_j 
\]
linear $\de \e$, $\de$ and $\e$-polynomials, respectively.
The
term, which according to the order defined previously is ordered highly among the
terms of $A$, is called the {\it maximum} term of $A$. By $d(A)$ we denote the minimum
delay of $A$, in other words
 $d(A)=\min(\min({\bf i}),\min({\bf j})), ({\bf i,j}) \in {\bf I}_n \times {\bf J}_m$,
$n=0,\ldots,\nu$, $m=0,\ldots,\mu$.
The largest of the numbers
 $deg({\bf i}+{\bf j})$, is called the {\it degree} of $A$, denoted also by $deg(A)$. 
 The equality of two $\de \e$-polynomials is defined with the same manner with the $\de$-polynomials.
 
 \begin{exa}
For the polynomial: $ A=5\de_0+6\de_1+\de_1^2-\de_2\de_3^2+2\e_1-3\e_2+\e_1^2-2\de_1\e_1+3\de_1^2\e_1\e_2^3$
we have
${\bf I}_0=\{e\}$, 
${\bf I}_1=\{0,1\}$,
${\bf I}_2=\{(1,1)\}$, ${\bf I}_3=\{(2,3,3)\}$, ${\bf J}_0=\{0\}$, ${\bf J}_1=\{1,2\}$,
${\bf J}_2=\{(1,1)\}$, ${\bf J}_3=\{(1,2,2)\}$
and
$c_{0e}=5$, $c_{1e}=6$, $c_{(1,1)e}=1$, $c_{(2,3,3)e}=-1$, $c_{e1}=2$, $c_{e2}=-3$, $c_{e(1,1)}=1$,
$c_{11}=-2$, $c_{(1,1)(1,2,2)}=3$, $c_{\bf ij}=0$ for all the other cases. 
$d(A)=0$ and $\deg (A)=\deg(\de_1^2)+\deg (\e_1\e_2^3)=6$.
When this polynomial acts to a pair of sequences, we get:
$A[y(t),u(t)]=5y(t)+6y(t-1)+y^2(t-1)$$-y(t-2)y(t-3)^2+2u(t-1)-$
$3u(t-2)+u^2(t-1)$$-2y(t-1)u(t-1)+3y^2(t-1)u(t-1)u^3(t-2)$.
\end{exa}
Now, working analogously with the case of $\de$-polynomials we are going to define the two internal operations, the dot and the star-product for $\de\e$-polynomials. Let
\[ A=\sum_{n=0}^{\nu_1} \sum_{m=0}^{\mu_1} \sum_{({\bf i}_1,{\bf j}_1) \in {\bf I}_{n,1} \times {\bf J}_{m,1}} a_{\bf i_1j_1}\de_{{\bf i}_1}\e_{{\bf j}_1} \quad \mbox{and} \quad 
 B=\sum_{n=0}^{\nu_2} \sum_{m=0}^{\mu_2} \sum_{({\bf i}_2,{\bf j}_2) \in {\bf I}_{n,2} \times {\bf J}_{m,2}} b_{\bf i_2j_2}\de_{{\bf i}_2}\e_{{\bf j}_2}\]
be two $\de\e$-polynomials. Their dot-product, denoted by $A \cdot B$, is the $\de \e$-polynomial, defined as:
\[ A \cdot B=\sum_{n=0}^{\max(\nu_1,\nu_2)}\sum_{m=0}^{\max(\mu_1,\mu_2)}\sum_{({\bf i},{\bf j}) \in {\bf I}_n \times {\bf J}_m }c_{\bf i j}\de_{\bf i}\e_{\bf j}\]
with $\de_{\bf i}=\de_{{\bf i}_1} \cdot \de_{{\bf i}_2}$,
$\e_{\bf j}=\e_{{\bf j}_1} \cdot \e_{{\bf j}_2}$
and
$c_{\bf ij}=a_{\bf i_1j_1} \cdot b_{\bf i_2j_2}$. This dot-product corresponds also to the usual product among polynomials of many variables.
The star-product is also defined here as a composition of maps. Indeed, let $B,C,A$ be $\de \e$-polynomials. We define the maps:
$ B: F\times F \to F, [w(t),v(t)] \to y(t)=B[w(t),v(t)] $,
$ C: F\times F \to F , [w(t),v(t)] \to u(t)=C[w(t),v(t)] $,
$ A: F\times F \to F , [y(t),u(t)] \to z(t)=A[y(t),u(t)] $
The $\de\e$-polynomial which corresponds to the composition:
\[ A \circ [B,C]: F\times F \to F \quad , \quad [w(t),v(t)] \to z(t) \]
is called the star-product of the polynomials $A,B,C$ and it is denoted by $A \ast [B,C] $.
The above definition can be either reduced to some special cases, when $B$ is only a $\de$-polynomial and $C$ an $\e$-one, or extended to more general situations, so that substitutions of a pair of $\de \e$-polynomials into a pair of $\de\e$-polynomials, to be described. Details can be found at \cite{kn:lappas1}.
A question, which naturally arises, is how can we calculate the star-product of $\de\e$-polynomials. There are certain formulas, \cite{kn:bibo},\cite{kn:lappas1}. The most general one is presenting in \cite{kn:lappas1}, where we can calculate the star-product between pairs of $\de\e$-polynomials, by means of tensor-products. Nevertheless, for the problem we study in the current paper we do not need this " huge " relation. The formulas, which are sufficient to cover all the cases we meet, are the following:

\begin{pro}\label{star}
Let $A$ and $B$ be two $\de\e$-polynomials, with $ A=\sum_{n=0}^\nu\sum_{m=0}^\mu\sum_{({\bf i,j}) \in {\bf I}_n \times {\bf J}_m}c_{\bf ij}\de_{\bf i}\e_{\bf j} $, then:
\par
\noindent
(a) $ \de_{i} \ast A=\sum_{n=0}^\nu\sum_{m=0}^\mu\sum_{({\bf i,j}) \in {\bf I}_n \times {\bf J}_m}c_{\bf ij}\de_{{\bf i} \dot{+} i}\e_{{\bf j} \dot{+} i}$, $\quad$ $i \in {\bf Z}^+$.
\par
\noindent
(b) $\de_{\bf i} \ast A=(\de_{i_1}\ast A)\cdot(\de_{i_2}\ast A)\cdots (\de_{i_k}\ast A)$, 
${\bf i}=(i_1,i_2,\ldots,i_k)$.
\par
\noindent
(c) $\de_{\bf i}\e_{\bf j} \ast [A,B]=(\de_{\bf i} \ast A) \cdot (\e_{\bf j} \ast B) $,
 where by 
" $\cdot$ " we denote the dot-product.
\par
\noindent
(d) Let $\de_{\bf i}\e_{\bf j}$ be a $\de$-operator, with ${\bf
i}=(i_1,i_2,\ldots,i_{\varphi})$, ${\bf j}=(j_1,j_2,\ldots,j_{\lambda})$ and
$L=\sum_{i \in I_1}r_i \de_i $, $M=\sum_{j \in J_1} m_j \e_j$ linear $\de$ and
$\e$-polynomials, $I_1,J_1 \subset Z$ finite sets of multindices, then
\[ \de_{\bf i}\e_{\bf j} \ast [L,M]=\sum_{{\bf s} \in I_1^{\varphi}}
\sum_{{\bf h} \in J_1^{\lambda}}r_{s_1}r_{s_2}\cdots r_{s_\varphi} m_{h_1} m_{h_2}
\cdots
 m_{h_\lambda}\de_{{\bf s}{+}{\bf i}}\e_{{\bf h}{+}{\bf j}} \]
 where ${\bf s}=(s_1,s_2,\ldots, s_{\varphi})$, ${\bf h}=(h_1,h_2, \ldots, h_{\lambda})$,
 $I_1^{\varphi}, J_1^{\lambda}$ the cartesian products
 \newline
 $I_1^{\varphi}=\underbrace{I_1\times I_1 \times \cdots \times I_1}_{\varphi - times}$,
 $J_1^{\lambda}=\underbrace{J_1\times J_1 \times \cdots \times J_1}_{\lambda - times}$, and
 ${\bf s}{+}{\bf i}$, ${\bf h}{+}{\bf j}$ the common vector addition.
 \end{pro}
{\bf Proof:} The proofs of (a), (b) and (c) are straightforward. We use these properties in order to show (d). Successively we have:
\[ \de_{\bf i}\e_{\bf j} \ast [L,M]=(\de_{\bf i} \ast L) \cdot (\e_{\bf j} \ast M)=
(\de_{i_1} \ast L)\cdots (\de_{i_\varphi} \ast L)\cdot (\e_{j_1} \ast M) \cdots
(\e_{j_\lambda} \ast M)=\]
\[ =\left(\sum_{i\in  I_1}r_i\de_{i+i_1}\right) \cdot \left(\sum_{i\in  I_1}r_i\de_{i+i_2}\right) \cdots \left(\sum_{i\in {\bf I}_1}r_i\de_{i+i_\varphi}\right)
\cdot\left(\sum_{j \in J_1} m_j \e_{j+j_1}\right)
\cdot\]
\[\cdot\left(\sum_{j \in J_1} m_j \e_{j+j_2}\right)
\cdots\left(\sum_{j \in J_1} m_j \e_{j+j_\lambda}\right)\]
The given formula is nothing else than a reinstall of the above products in a compact form.
\begin{exa}
 (1) $ \de_1\de_2 \e_2^2 \ast [2\de_0+\de_1,2\e_0-\e_1]=$
 $[\de_1\de_2 \ast (2\de_0+\de_1) ] \cdot [\e_2^2\ast(2\e_0-\e_1)]=$
 $(2\de_1+\de_2) \cdot (2\de_2+\de_3) $
 $\cdot (2\e_2-\e_3) \cdot (2\e_2-\e_3)=$
 $16\e_2^2\de_1\de_2 - 16\e_2\e_3\de_1\de_2 +$
$ 4\e_3^2\de_1\de_2 + 8\e_2^2\de_2^2 - 8\e_2\e_3\de_2^2 +$ $2\e_3^2\de_2^2 +
8\e_2^2\de_1\de_3 - $ $ 8\e_2\e_3\de_1\de_3 + 2\e_3^2\de_1\de_3 +$ $ 4\e_2^2\de_2\de_3
- 4\e_2\e_3\de_2\de_3 + \e_3^2\de_2\de_3$
\par
(2) $(2 \de_1\e_2+\de_0^2\e_1\e_2)$ $ \ast [\de_1-2\de_2, \e_0^2\e_1]$ $=2 \de_1\e_2
\ast [\de_1-2\de_2, \e_0^2\e_1]$ $+ \de_0^2\e_1\e_2 \ast [\de_1-2\de_2, \e_0^2\e_1]$
$=[2 \de_1 \ast (\de_1-2\de_2)] \cdot [\e_2 \ast \e_0^2\e_1]$ $+[\de_0^2 \ast
(\de_1-2\de_2)] \cdot [\e_1\e_2 \ast  \e_0^2\e_1]$ $=2(\de_2-2\de_3) \cdot
(\e_2^2\e_3)+$ $(\de_1-2\de_2)^2 \cdot \e_1^2\e_2 \cdot \e_2^2\e_3=$
$2\de_2\e_2^2\e_3$$-4\de_3\e_2^2\e_3$$+\de_1^2\e_1^2\e_2^2\e_3$$-4\de_1\de_2\e_1^2\e_2^2\e_3$$+4\de_2^2\e_1^2\e_2^2\e_3$

\end{exa}

\section{The Formal Linear Factors}
In this section we present the main tool that is used throughout this paper. Actually, it is about detecting linear polynomials which can be considered as " factors " of a given $\de\e$-polynomial with respect to the star-product and the dot-product. Specifically, let $A$ be a concrete $\de\e$-polynomial we write it as follows:
\begin{equation}\label{compact}A=\sum_{\bf \theta} \de_0^{\kappa_{\bf \theta}}\e_0^{\sigma_{\bf \theta}}\cdot (c_{\bf \theta}\de_{{\bf i}_{\bf \theta}}\e_{{\bf j}_{\theta}} \ast [L_{\bf \theta},M_{\bf \theta}])+{\bf R}_{\de,h}+{\bf R}_{\e,h}+R_{\de\e}
\end{equation}
where
$L_{\bf \theta},M_{\bf \theta}$ are linear $\de$ and $\e$-polynomials, correspondingly,
(a detailed explanation of the above expression will be gived later), including the linear $\de$ and $\e$-parts of $A$ and $R_{\de,h},R_{\e,h}$,$R_{\de\e}$ the remainders.
The novelty is that some coefficients of those polynomials are not concrete numbers but undetermined parameters, which can take arbitrary values. We call these coefficients {\it formal coefficients} and the linear polynomials, {\it formal linear factors} of $A$.
The cornerstone of the current paper is relied on the fact that by giving to these parameters proper values, we can take certain expressions of $A$, appropriate for 
solving the model matching problem. 
\par\noindent
We elucidate here that we get the expression (\ref{compact}) and the sets of formal linear factors by means of a symbolic algorithm. 
This algorithm works as follows. First it decomposes $A$ to its $\de$, $\e$ and $\de\e$-parts. For each such part it finds the corresponding maximum term and then it constructs formal linear polynomials $L_\theta,M_\theta$ and an operator $\de_{{\bf i}_\theta}\e_{{\bf j}_\theta}$, so that the operation $\de_{{\bf i}_\theta}\e_{{\bf j}_\theta} \ast [L_\theta,M_\theta]$ eliminates it. Working recursively we eliminate all the other terms but a remainder ${\bf R}$. This remainder is factorized with respect to the dot-product as $R=\de_0^k\e_0^\sigma \cdot \tilde{R}$ and then we repeat the whole procedure, working with $\tilde{R}$. To establish certain definitions to a more rigorous way we need first to present the said algorithm.

\small {\sf
\par
\vskip 10 pt \noindent \underline{\bf The Formal Linear Factors (FLF) - Algorithm.}
\par
\vskip 10 pt {\bf Input:} A $\de \e$-polynomial $A$ with no constant terms.
\par
\vskip 5 pt {\bf Output:} The quantity  $h^*$, the sets ${\cal L},{\cal M},{\cal C},{\cal O}$ and the polynomials $R_{\de,\nu},R_{\e,\nu}$, $\nu=0,1,\ldots,h^*$, $R_{\de\e,h^*}$.
\vskip 8 pt
\noindent
{\bf STEP 0:} Set $h=0$, ${\cal L}=\{\}$, ${\cal M}=\{ \}$, ${\cal C}=\{\}$, ${\cal O}=\{ \}$.
\par
\noindent
{\bf STEP 1:} Decompose $A$ as follows: 
$A=A_{\de}+A_{\e}+A_{\de\e}$, 
where $A_{\de}, A_{\e},
A_{\de\e}$ are the pure $\de$-part, $\e$-parts and $\de\e$- parts of
$A$, respectively.
\par
\noindent {\bf STEP 2:} Set $u=1,\omega=0$. {\bf REPEAT} the following substeps {\bf UNTIL} $A_{\de}=0$ or $A_{\de}$ is a polynomial of the variable $\de_0$ only.
\begin{quote}
\par
\noindent {\bf SUBSTEP 2a:} {\bf CALL SUBROUTINE REMAINDER} with inputs $A_{\de}, {\cal L}, {\cal M}, {\cal C}, {\cal O}$, $h,u,\omega$ and outputs $R,{\cal L}, {\cal M}, {\cal C}, {\cal O}$. 
\par
\noindent {\bf SUBSTEP 2b:} Find the maximum integer $\theta \ge 0$ such that $R=\de_0^\theta \cdot \tilde{R}$ and $\tilde{R}$ is a $\de$-polynomial, not containing a constant term.
\par
\noindent {\bf SUBSTEP 2c:} Set $A_{\de}=\tilde{R}$, $u=u+1$.

\end{quote}
{\bf END OF THE REPEAT} 
\par
\noindent
{\bf STEP 3:} Rename the last value $A_{\de}$ as $R_{\de,h}$.
\par
\noindent {\bf STEP 4:} Set $u=0,\omega=1$. {\bf REPEAT} the following substeps {\bf UNTIL} $A_{\e}=0$ or $A_{\e}$ is a polynomial of the variable $\e_0$ only.
\begin{quote}
\par
\noindent {\bf SUBSTEP 4a:}  {\bf CALL SUBROUTINE REMAINDER} with input $A_{\e}, {\cal L}, {\cal M}, {\cal C}, {\cal O}$,$h,u,\omega$ and outputs $R,{\cal L}, {\cal M}, {\cal C}, {\cal O}$.
\par
\noindent {\bf SUBSTEP 4b:} Find the maximum integer $\rho \ge 0$ such that $R=\e_0^\rho \cdot \tilde{R}$ and $\tilde{R}$ is an $\e$-polynomial, not containing a constant term.
\par
\noindent {\bf SUBSTEP 4c:} Set $A_{\e}=\tilde{R}$, $\omega=\omega+1$.

\end{quote}
{\bf END OF THE REPEAT}
\par
\noindent
{\bf STEP 5:} Rename the last value $A_{\e}$ as $R_{\e,h}$.
\par
\noindent
{\bf STEP 6:} {\bf IF} $A_{\de\e}=0$ or $A_{\de\e}$ is a polynomial of the variables
$\de_0$ and $\e_0$ only {\bf THEN} set $R_{\de\e,h}=A_{\de\e}$, $h^*=h$
 and goto the {\bf Output} {\bf ELSE} proceed.
\par
\noindent
{\bf STEP 7:} Set $u=0,\omega=0$, {\bf CALL SUBROUTINE REMAINDER} with input $A_{\de\e}, {\cal L}, {\cal M}$, ${\cal C}$, ${\cal O}$, $h,u,\omega$ and outputs $R,{\cal L}, {\cal M}, {\cal C}, {\cal O}$.
\par
\noindent
{\bf STEP 8:} Find the maximum integers $\kappa,\sigma \ge 0$ such that $R=\de_0^\kappa\e_0^\sigma\cdot \tilde{R}$ and $\tilde{R}$ is a $\de\e$-polynomial, not containing a constant term.
\par
\noindent
{\bf STEP 9:} Set $A=\tilde{R}$, $h=h+1$, goto step 1.

\par
\vskip 15 pt \noindent \underline{\bf SUBROUTINE REMAINDER.}
\par
\vskip 15 pt {\bf Input:} A $\de \e$-polynomial $A$,  the sets ${\cal L},{\cal M}, {\cal C}, {\cal O}$ and the counters $h,u,\omega$.
\par
\vskip 5 pt {\bf Output:} The polynomial  $R$. The sets ${\cal L}, {\cal M}, {\cal C}, {\cal O}$.
\vskip 10 pt
\begin{quote}
\par
\noindent
{\bf REM-STEP 0:} Set $\lambda=0$, $R=A$.
\par
\noindent
{\bf REM-STEP 1:} We decompose $R$ as follows:
$R=R_{l}+R_{nl}$
where $R_{l}$,$R_{nl}$ are the linear and the non-linear parts of $R$, respectively.
\par
\noindent
{\bf REM-STEP 2:} {\bf IF} $R_l \ne 0$  {\bf THEN}
{\bf IF } $u=1$ {\bf THEN }  set
$L_{0,u,h}=R_{l}$, $c_{0,u,\omega,h}=1$, ${\cal O}={\cal O} \cup \{\de_0\}$, ${\cal C}={\cal C} \cup \{1\}$, ${\cal L}={\cal L} \cup \{L_{0,u,h}\}$ {\bf ELSE} set $M_{0,\omega,h}=R_{l}$, 
$c_{0,u,\omega,h}=1$, ${\cal O}={\cal O} \cup \{\de_0\}$, ${\cal C}={\cal C} \cup \{1\}$,
${\cal M}={\cal M} \cup \{M_{0,\omega,h}\}$ {\bf ELSE} proceed.
\par
\noindent {\bf REM-STEP 3:}  Set $R=R_{nl}$ {\bf REPEAT}  the following REM-substeps {\bf UNTIL} $R_{\lambda}=0$ or all the terms of $R_{\lambda}$ become zero terms (i.e. terms of the form
$\de_0\de_{i_2}\cdots \de_{i_n}$$\e_0\e_{j_1}\cdots\e_{j_m}$,
$0\le i_2 \le \cdots \le i_n$, $0 \le j_2 \le \cdots j_m$).

\begin{quote}
\par
\noindent{\bf REM-Substep 3a:} $\lambda=\lambda+1$, $R_{\lambda-1}=R$
\par
\noindent{\bf REM-Substep 3b:} Let $c_{\lambda,u,\omega,h}$$\de_{{\bf i}_{\lambda,u,\omega,h}}$$\e_{{\bf j}_{\lambda,u,\omega,h}}=$
$c_{\lambda,u,\omega,h} \de_{i_1}\de_{i_2} \cdots \de_{i_n} \e_{j_1}\e_{j_2} \cdots
\e_{j_m} $ be the {\it maximum} non-zero term of $R_{\lambda-1}$, $0 \le i_1 \le i_2 \le  \cdots
\le i_n$, $0 \le j_1 \le j_2 \le  \cdots \le j_m$ and $c_{\lambda,u,\omega,h}$
its coefficient. (Actually, in the first iteration of the algorithm, $c_{\lambda,u,\omega,h}$ is always a real number. It then becomes a function of the
unknown parameters $w_{i,j,k,h},s_{i,j,k,h}$.)
\par
\noindent{\bf REM-Substep 3c:} {\bf IF} $\de_{{\bf i}_{\lambda,u,\omega,h}}\ne \de_e$ {\bf THEN} we form the formal linear $\de$-polynomial
\[ L_{\lambda,u,h}=w_{\lambda,u,h,0}\de_0+w_{\lambda,u,h,1}\de_1+
w_{\lambda,u,h,2}\de_2 +\cdots +\de_{i_1} \]
where $w_{\lambda,u,h,0},w_{\lambda,u,h,1},\ldots$ 
are unknown parameters, taking values in
${\bf R}$ and $i_1$ is the minimum delay of the $\de$-part of the maximum term {\bf ELSE} $L_{\lambda,u,h}=\de_e$
\par
\noindent{\bf REM-Substep 3d:} {\bf IF} $\e_{{\bf j}_{\lambda,u,\omega,h}}\ne \e_e$ {\bf THEN} 
we form the formal linear $\e$-polynomial
\[ M_{\lambda,\omega,h}=s_{\lambda,\omega,h,0}\e_0+s_{\lambda,\omega,h,1}\e_1+s_{\lambda,\omega,h,2}\e_2 
+\cdots +\e_{j_1} \]
where $s_{\lambda,\omega,h,0},s_{\lambda,\omega,h,1},\ldots$ are unknown parameters that take values in
${\bf R}$ and $j_1$ is the minimum delay of the $\e$-part of the maximum term {\bf ELSE}
$M_{\lambda,\omega,h}=\e_e$
\par
\noindent{\bf REM-Substep 3e:} We execute the operation:
\[ R_{\lambda}=R_{\lambda-1}-c_{\lambda,u,\omega,h}\de_{0}\de_{i_2-i_1} \cdots \de_{i_n-i_1}
\e_{0}\e_{j_2-j_1} \cdots \e_{j_m-j_1}
 \ast [L_{\lambda,u,h}, M_{\lambda,\omega,h}]=\]
 \[=R_{\lambda-1}-c_{\lambda,u,\omega,h}\de_{{\bf i}_{\lambda,u,\omega,h}}\e_{{\bf j}_{\lambda,u,\omega,h}} \ast [L_{\lambda,u,h}, M_{\lambda,\omega,h}]\]
 \par
\noindent{\bf REM-Substep 3f:} For $L_{\lambda,u,h}\ne \de_e$, $M_{\lambda,\omega,h}\ne \e_e$ we set ${\cal L}={\cal L}\cup \{L_{\lambda,u,h}\}$,
 ${\cal M}={\cal M}\cup \{M_{\lambda,\omega,h}\}$, ${\cal C}={\cal C} \cup \{ c_{\lambda,u,\omega,h}\}$,
 ${\cal O}={\cal O}\cup \{\de_{{\bf i}_{\lambda,u,\omega,h}}\e_{{\bf j}_{\lambda,u,\omega,h}}\}$
\end{quote}
{\bf END OF THE REPEAT}
\par
\noindent {\bf REM-STEP 4:} Rename the last value of $R_{\lambda}$ as $R$ and go to the {\bf Output}. 
\end{quote}}
\noindent
\normalsize
The following definitions will play a crucial role in the context of the model
matching problem, below. First of all
the sets
${\cal L}$ and ${\cal M}$, which are formed from all the linear $\de$ and
$\e$-polynomials that appear in the above algorithm, are called 'the sets of the
Formal $\de$-Linear Factors' and 'the Formal $\e$-Linear Factors of $A$', respectively.
We must note here that the use of the word Formal is a little bit excessive, since all the elements of the sets ${\cal L}$ and ${\cal M}$ are not formal linear polynomials. The linear $\de$ and $\e$-parts of $A$, usually denoted by $A_{\de,l},A_{\e,l}$, are a first example. The linear polynomial $R_l$, appeared in REM-STEP 1 is another. By ${\cal L}^*$ and ${\cal M}^*$ we name the sets
${\cal L}^*={\cal L}-\{A_{\de,l}\}$, ${\cal M}^*={\cal M}-\{A_{\e,l}\}$, ${\cal L}_{\de}$ is the set of the linear factors of the $\de$-part of $A$, in other words ${\cal L}_{\de}=\{ L_{\lambda,u,0}\}\subset {\cal L}$
(it includes $A_{\de,l}$), $\overline{\cal L}$, $\overline{\cal M}$ are the sets of the $\de$ or $\e$-linear factors of $A$ which do not participate to any cross-product, i.e. 
$\overline{\cal L}=\{ L_{\lambda,u,h}, \quad u \ne 0\}$,
$\overline{\cal M}=\{ M_{\lambda,\omega,h}, \quad \omega \ne 0\}$, finally
$\overline{\cal L}^*=\overline{\cal L}-\{A_{\de,l}\}$,
$\overline{\cal M}^*=\overline{\cal M}-\{A_{\e,l}\}$. The coefficient $c_{\lambda,u,\omega,h}$ is called the {\it companion coefficient} of the polynomial $L_{\lambda,u,h}$ or of the polynomial $M_{\lambda,\omega,h}$ and it is denoted by $ccoef(L_{\lambda,u,h})$ or $ccoef(M_{\lambda,\omega,h})$. By convention,
$ccoef(L_{0,u,h})=ccoef(M_{0,\omega,h})=1$ for any $u$ and $h$. The set of all companion coefficients is denoted by ${\cal C}$.
The polynomials $R_{\de,\nu},R_{\e,\nu}$, $\nu=0,1,\ldots,h^*$, $R_{\de\e,h^*}$ are called the $\de$,$\e$ and $\de\e$-remainders, correspondingly and the set ${\cal O}$ the set of the formal $\de\e$-operators. Let ${\cal A} \subseteq {\cal L} \cup {\cal M}$, be a set consisting from some $\de$ or $\e$-linear factors of $A$. By ${\cal C}^{\cal A}\subseteq {\cal C}$ we denote the set of the companion coefficients of the elements of ${\cal A}$. That is
${\cal C}^{\cal A}=\{c_{\lambda,u,\omega,h} \in {\cal C}: c_{\lambda,u,\omega,h} $ is the companion coefficient of an element of ${\cal A}\}$.
\noindent
The next theorems describe the behaviour of the FLF-Algorithm.
\begin{theo}
The FLF-Algorithm terminates after a finite number of steps.
\end{theo}
{\bf Proof:} A glance at the algorithm indicates that essentially, it consists from three repeat-loops, appeared at the steps 2,4 and the REM-step 3 and a loop defined between the steps 1 and 9.
The termination of those loops ensure the termination of the algorithm. Let $A$ be an input of the Remainder-Subroutine. As we pass from the REM-substep 3e we get:
$ R_{new}=R_{previous}-c_{\lambda,u,\omega,h}\de_{{\bf i}_{\lambda,u,h}}\e_{{\bf j}_{\lambda,\omega,h}}\ast[L_{\lambda,u,h},M_{\lambda,\omega,h}]$, where $c_{\lambda,u,\omega,h}$ is the current maximum non zero term. By means of proposition
\ref{star} (d), we get:
\[ \de_{{\bf i}_{\lambda,u,h}}\e_{{\bf j}_{\lambda,\omega,h}}\ast[L_{\lambda,u,h},M_{\lambda,\omega,h}]=\de_{0}\de_{i_2-i_1} \cdots \de_{i_n-i_1}
\e_{0}\e_{j_2-j_1}\cdots \e_{j_m-j_1}\ast[L_{\lambda,u,h},M_{\lambda,\omega,h}]=\]
\[=\de_{0}\de_{i_2-i_1} \cdots \de_{i_n-i_1}
\e_{0}\e_{j_2-j_1}\cdots \e_{j_m-j_1}\ast [w_{\lambda,u,h,0}\de_0+\cdots +\de_{i_1},s_{\lambda,\omega,h,0}\e_0+\cdots +\e_{j_1}]=\]
\[=\de_{i_1}\de_{i_2} \cdots \de_{i_n}\e_{j_1}\e_{j_2}\cdots \e_{j_m} +\mbox{other terms}\]
Therefore, the above substraction will eliminate this $\de_{{\bf i}_{\lambda,u,h}}\e_{{\bf j}_{\lambda,\omega,h}}$ term. The construction of the linear polynomial $L_{\lambda,u,h},M_{\lambda,\omega,h}$ guarantees that all the terms of $R_{new}$ will be ordered lower than $\de_{\bf i}\e_{\bf j}$, thus, by repeating the same procedure for those terms, we shall gradually eliminate all the non zero terms of $A$ and the REMAINDER 
SUBROUTINE will terminate. Let $R$ be the output of this subroutine. It will be consist of terms
of the form: $\de_0\de_{\bf k}\e_o\e_{\bf h}$. Factorizing it by powers of $\de_0\e_0$ (steps 2b,4b,8) we decrease the degree of $R$ at least by two.
This nesting ensures the termination of all the loops and consequently the termination of the algorithm.

\begin{theo}\label{typos}
Let $A$  be a given $\de\e$-polynomial and $A_{l,\de}$, $A_{l,\e}$ its linear $\de$ and $\e$-parts. 
Let ${\cal L,M,C,O}$, $h^*$ and $R_{\de\e,h^*}$, $R_{\de,h},R_{\e,h}$, $h=0,\ldots,h^*$ be the outputs of the FLF-Algorithm, after its application to a given $\de\e$-polynomial $A$. Then,
\begin{equation}\label{FLFMAIN}
A=\sum_{h=0}^{h^*}\sum_{\omega=0}^{\omega_h^*}\sum_{u=0}^{u_h^*}\sum_{\lambda=0}^{v_{u,\omega,h}} \de_0^{\kappa_{u,h}}\e_0^{\sigma_{\omega,h}} \cdot (c_{\lambda,u,\omega,h}\de_{{\bf i}_{\lambda,u,\omega,h}}\e_{{\bf j}_{\lambda,u,\omega,h}}\ast [L_{\lambda,u,h},M_{\lambda,\omega,h}])+{\bf R}_{\de,h}+{\bf R}_{\e,h}+R_{\de\e,h^*}
\end{equation}
where $\omega_h^*,u_h^*,v_{u,\omega,h},\kappa_{u,h},\sigma_{\omega,h}$ are positive integers, depending from the values of $h,u$ and $\omega$, $c_{\lambda,u,\omega,h} \in {\cal C}$, $\de_{{\bf i}_{\lambda,u,\omega,h}}\e_{{\bf j}_{\lambda,u,\omega,h}} \in {\cal O}$,
$L_{\lambda,u,h} \in {\cal L}$, $M_{\lambda,\omega,h} \in {\cal M}$ and the following " border " conditions hold: $\kappa_{1,h}=0$, $\kappa_{0,0}=0$, $\sigma_{1,h}=0$, $\sigma_{0,0}=0$, $c_{0,0,\omega,h}=1$, $c_{0,u,0,h}=1$, $c_{0,0,0,h}=0 $, $c_{\lambda,u,\omega,h}=0$ for $u+\omega >1$, $\de_{{\bf i}_{0,u,0,h}}=\de_0 $, $\de_{{\bf i}_{\lambda,0,\omega,h}}=\de_e$,
$\e_{{\bf j}_{0,0,\omega,h}}=\e_0 $, $ \e_{{\bf j}_{\lambda,0,\omega,h}}=\e_e$,
for $u=1,\ldots,u_h^*$, $\omega=1,\ldots,\omega_h^*$, $\lambda=0,\ldots,v_{u,\omega,h}$, $h=0,\ldots,h^*-1$ and
$L_{0,1,0}=A_{\de,l}$, $M_{0,1,0}=A_{\e,l}$. Furthermore, this expression is feasible for any $\de\e$-polynomial $A$ and unique, under the assumption that the parameters $w_{\lambda,u,h,k}$, $s_{\lambda,\omega,h,k}$ are considered as constants.
\end{theo}
{\bf Proof:} Using reverse recursive relations we can write the evolution of the FLF-Algorithm  in a condensed form. Let us denote by $H_h$ the outcome of the loop, defined by steps 1 and 9, at any specific iteration $h$, $\Delta_{u,h}$ the outcome of the  repeat procedure in step 2, and $E_{\omega,h}$ the outcome at the repeat procedure in step 4. It can be easily proved that the next equalities hold:
\[ H_h=\Delta_{1,h}+E_{1,h}+\sum_{\lambda=1}^{\nu_h}c_{\lambda,0,0,h}\de_{{\bf i}_{\lambda,0,0,h}}\e_{{\bf j}_{\lambda,0,0,h}} \ast [L_{\lambda,0,h},M_{\lambda,0,h}] +\de_0^{\kappa_{u,h+1}}\e_0^{\sigma_{\omega,h+1}} \cdot H_{h+1}\quad,\quad 
h=0,1, \ldots,h^*-1\]
\[ \Delta_{u,h}=L_{0,u,h}+\sum_{\lambda=1}^{g_{u,h}}c_{\lambda,u,0,h}\de_{{\bf i}_{\lambda},u,h}\ast L_{\lambda,u,h}+\de_0^{k_{u,h}}\cdot \Delta_{u+1,h}\quad,  \quad u=1,\ldots, u_h^* \]
\[ E_{\omega,h}=M_{0,\omega,h}+\sum_{\lambda=1}^{\gamma_{\omega,h}}c_{\lambda,0,\omega,h}\e_{{\bf j}_{\lambda},\omega,h}\ast M_{\lambda,\omega,h}+\e_0^{\sigma_{\omega,h}}\cdot E_{\omega+1,h} \quad, \quad \omega=1,\ldots, \omega_h^* \]
with terminal conditions $H_{h^*}=R_{\de\e,h^*}$, $\Delta_{u^*_h+1,h}=R_{\de,h}$, $E_{\omega^*_h+1,h}=R_{\e,h}$, $h=0,1,\ldots,h^*$, $\nu_h$, $g_{u,h}$, $\gamma_{\omega,h}$ are positive integers, depending from the quantities $h,u,\omega$.
Obviously $A=H_0$. Expanding this relation in details we take:
\[ A=L_{0,1,0}+\sum_{\lambda=1}^{g_{1,0}}c_{\lambda,1,0,0}\de_{{\bf i}_{\lambda,1,0}}\ast L_{\lambda,1,0}+\de_0^{k_{1,0}}\cdot (L_{0,2,0}+\sum_{\lambda=1}^{g_{2,0}}c_{\lambda,2,0,0}\de_{{\bf i}_{\lambda,2,0}}\ast L_{\lambda,2,0}+\de_0^{k_{2,0}}\cdot ( \]
\[ \cdots \de_0^{k_{u^*_0-1,0}}\cdot (L_{0,u^*_0,0}+\sum_{\lambda=1}^{g_{u^*_0,0}}c_{\lambda,u^*_0,0,0}\de_{{\bf i}_{\lambda},u^*_0,0}\ast L_{\lambda,u^*_0,0}+ R_{\de,0} )\cdots\left. \right)+M_{0,1,0}+\sum_{\lambda=1}^{\gamma_{1,0}}c_{\lambda,0,1,0}\]
\[\e_{{\bf j}_{\lambda,1,0}}\ast M_{\lambda,1,0}+\e_0^{\sigma_{1,0}}\cdot(M_{0,2,0}+\sum_{\lambda=1}^{\gamma_{2,0}}c_{\lambda,0,2,0}\e_{{\bf j}_{\lambda,2,0}}\ast M_{\lambda,2,0}+\e_0^{\sigma_{2,0}}\cdot(
\cdots  \e_0^{\sigma_{\omega^*_0-1,0}}\cdot(  M_{0,\omega^*_0,0}+\]

\[+\sum_{\lambda=1}^{\gamma_{\omega^*_0,0}}c_{\lambda,0,\omega^*_0,0}\e_{{\bf j}_{\lambda,\omega^*_0,0}}\ast M_{\lambda,\omega^*_0,0}+R_{\e,0} )) \cdots )+
\sum_{\lambda=1}^{\nu_0}c_{\lambda,0,0,0}\de_{{\bf i}_{\lambda,0,0,0}}\e_{{\bf j}_{\lambda,0,0,0}} \ast [L_{\lambda,0,0},M_{\lambda,0,0}] +\]

\[+\de_0^{\kappa_{0,1}}\e_0^{\sigma_{0,1}} \cdot ( L_{0,1,1}+\sum_{\lambda=1}^{g_{1,1}}c_{\lambda,1,0,1}\de_{{\bf i}_{\lambda,1,1}}\ast L_{\lambda,1,1}+\de_0^{k_{1,1}}\cdot (L_{0,2,1}+\sum_{\lambda=1}^{g_{2,1}}c_{\lambda,2,0,1}\de_{{\bf i}_{\lambda,2,1}}\ast L_{\lambda,2,1}+\de_0^{k_{2,1}}\cdot  ( \]

\[ \cdots \de_0^{k_{u^*_1-1,1}}\cdot (L_{0,u^*_1,1}+\sum_{\lambda=1}^{g_{u^*_1,1}}c_{\lambda,u^*_1,0,1}\de_{{\bf i}_{\lambda,u^*_1,1}}\ast L_{\lambda,u^*_1,1}+ R_{\de,1}) )\cdots )+M_{0,1,1}+\sum_{\lambda=1}^{\gamma_{1,1}}c_{\lambda,0,1,1}\e_{{\bf j}_{\lambda,1,1}}\ast M_{\lambda,1,1}+\]
\[+\e_0^{\sigma_{1,1}}\cdot (M_{0,2,1}+\sum_{\lambda=1}^{\gamma_{2,1}}c_{\lambda,0,2,1}\e_{{\bf j}_{\lambda,2,1}}\ast M_{\lambda,2,1}+\e_0^{\sigma_{2,1}} (
\cdots \e_0^{\sigma_{\omega^*_1-1,1}}(  M_{0,\omega^*_1,1}+\sum_{\lambda=1}^{\gamma_{\omega^*_1,1}}c_{\lambda,0,\omega^*_1,1}\e_{{\bf j}_{\lambda},\omega^*_1,1}\ast M_{\lambda,\omega^*_1,1}+\]

\[+R_{\e,1} )) \cdots )+\sum_{\lambda=1}^{\nu_1}c_{\lambda,0,0,1}\de_{{\bf i}_{\lambda,0,0,1}}\e_{{\bf j}_{\lambda,0,0,1}} \ast [L_{\lambda,0,1},M_{\lambda,0,1}] +\de_0^{\kappa_{0,2}}\e_0^{\sigma_{0,2}} (
 \cdots \cdots
 L_{0,1,h^*-1}+\sum_{\lambda=1}^{g_{1,h^*-1}}c_{\lambda,1,0,h^*-1}\]
 \[\de_{{\bf i}_{\lambda,1,h^*-1}}\ast L_{\lambda,1,h^*-1}+\de_0^{k_{1,h^*-1}} (L_{0,2,h^*-1}+\sum_{\lambda=1}^{g_{2,h^*-1}}c_{\lambda,2,0,h^*-1}\de_{{\bf i}_{\lambda,2,h^*-1}}\ast L_{\lambda,2,h^*-1}+\de_0^{k_{2,h^*-1}}( \]
  
\[ \cdots \de_0^{k_{u^*_{h^*-1}-1,h^*-1}} (L_{0,u^*_{h^*-1},h^*-1}+\sum_{\lambda=1}^{g_{u^*_{h^*-1},h^*-1}}c_{\lambda,u^*_{h^*-1},0,h^*-1}\de_{{\bf i}_{\lambda,u^*_{h^*-1},h^*-1}}\ast L_{\lambda,u^*_{h^*-1},h^*-1}+ R_{\de,h^*-1}))\cdots\]

\[)+M_{0,1,h^*-1}+\sum_{\lambda=1}^{\gamma_{1,h^*-1}}c_{\lambda,0,1,h^*-1}\e_{{\bf j}_{\lambda,1,h^*-1}}\ast M_{\lambda,1,h^*-1}+\e_0^{\sigma_{1,h^*-1}}(M_{0,2,h^*-1}+\sum_{\lambda=1}^{\gamma_{2,h^*-1}}c_{\lambda,0,2,h^*-1}\e_{{\bf j}_{\lambda,2,h^*-1}}\ast \]

\[\ast M_{\lambda,2,h^*-1}+\e_0^{\sigma_{2,h^*-1}}( \cdots \e_0^{\sigma_{\omega^*_{h^*-1}-1,h^*-1}}(  M_{0,\omega^*_{h^*-1},h^*-1}+\sum_{\lambda=1}^{\gamma_{\omega^*_{h^*-1},h^*-1}}c_{\lambda,0,\omega^*_{h^*-1},h^*-1}\e_{{\bf j}_{\lambda,\omega^*_{h^*-1},h^*-1}}\ast\]
\[ \ast M_{\lambda,\omega^*_{h^*-1},h^*-1}+R_{\e,h^*-1} )) \cdots )+\]
\begin{equation}\label{big}
 +\sum_{\lambda=1}^{\nu_h^*-1}c_{\lambda,0,0,h^*-1}\de_{{\bf i}_{\lambda,0,0,h^*-1}}\e_{{\bf j}_{\lambda,0,0,h^*-1}} \ast [L_{\lambda,0,h^*-1},M_{\lambda,0,h^*-1}] +R_{\de\e,h^*} ))) \cdots ))\cdots )
 \end{equation}
Rewriting (\ref{big}) in a proper compact manner and renaming the quantities $g_{u,h}$, $\gamma_{u,h}$ as $v_{u,\omega,h}$, we take the relation (\ref{FLFMAIN}) upon proving.
Now, let $A$ be an arbitrary $\de\e$-polynomial. If it is a linear one, of the form (\ref{l}), then it can be written like (\ref{FLFMAIN}) in a trivial way, by setting $L_{0,1,0}=A_{\de,l}$, $M_{0,1,0}=A_{\e,l}$, $c_{\lambda,u,\omega,h}=0$ 
$h,\lambda \ne 0$, $u,\omega \ne 1$ and putting the remainders equal to zero. If $A$ is nonlinear but contains only $\de_0\e_0$-terms, then we put $R_{\de\e,h^*}=A$ and $c_{\lambda,u,\omega,h}=0$. In all the other cases we have a maximum term and therefore the FLF-Algorithm can be applied. All the above supports the feasibility of the method. Let us now deal with the uniqueness question. 
We explain first that by the expression " are considered as 
constants " we mean that we treat the parameters like to be specific number.
Now, if $A$ has only zero terms or it is linear then the proof of the theorem is trivial. Let us suppose that $A$ contains non-zero terms and that, for the sake of proof, 
$\psi_a\de_{{\bf i}_a} \e_{{\bf j}_a}, \psi_a \in {\bf R}$ with
${\bf i}_a=(i_{1,a},i_{2,a},\ldots,i_{n,a})$,
${\bf j}_a=(j_{1,a},j_{2,a},\ldots,j_{m,a})$, is its maximum term. Let us further suppose, for the sake of simplicity, that the second higher ordered term of $A$ is  
$\psi_\beta\de_{{\bf i}_\beta} \e_{{\bf j}_\beta}, \psi_\beta \in {\bf R}$ with
${\bf i}_{\beta}=(i_{1,\beta},i_{2,\beta},\ldots,i_{n,\beta})$,
${\bf j}_{\beta}=(j_{1,\beta},j_{2,\beta},\ldots,j_{m,\beta})$ and $j_{m,\beta}=j_{m,a}-1$. As we pass from the REM-Step 3e for the first time we get:
$R_1=A-c_{1,0,0,0}\de_{{\bf i}_{1,0,0,0}}\e_{{\bf j}_{1,0,0,0}}$$\ast[L_{1,0,0},M_{1,0,0}]$ and thus $A=R_1+c_{1,0,0,0}\de_{{\bf i}_{1,0,0,0}}\e_{{\bf j}_{1,0,0,0}}$$\ast[L_{1,0,0},M_{1,0,0}]$. But 
$\de_{{\bf i}_{1,0,0,0}}\e_{{\bf j}_{1,0,0,0}}$
$=\de_0\de_{i_{2,a}-i_{1,a}}\cdots \de_{i_{n,a}-i_{1,a}}$
$\e_0\e_{j_{2,a}-j_{1,a}}\cdots \e_{j_{m,a}-j_{1,a}}$ and
$L_{1,0,0}=w_{1,0,0,0}\de_0+w_{1,0,0,1}\de_1+\cdots +\de_{i_{1,a}}$,
$M_{1,0,0}=s_{1,0,0,0}\e_0+s_{1,0,0,1}\e_1+\cdots +\e_{j_{1,a}}$. By means 
of proposition
(\ref{star}), (d), we can see that 
$c_{1,0,0,0}\de_{{\bf i}_{1,0,0,0}}\e_{{\bf j}_{1,0,0,0}}$$\ast[L_{1,0,0},M_{1,0,0}]$
$=c_{1,0,0,0}\de_{i_{1,a}}\de_{i_{2,a}}\cdots\de_{i_{n,a}}$
$\e_{j_{1,a}}\e_{j_{2,a}}\cdots\e_{j_{m,a}}$
$+c_{1,0,0,0}(s_{1,0,0,j_{m,a}-1})\de_{i_{1,a}}\de_{i,{2,a}}\cdots\de_{i_{n,a}}$
$\e_{j_{1,a}}\e_{j_{2,a}}\cdots\e_{j_{m,a}-1}$
$+$ other terms
$=c_{1,0,0,0}\de_{{\bf i}_a}\e_{{\bf j}_a}$
$+c_{1,0,0,0}(s_{1,0,0,j_{m,a}-1})\de_{{\bf i}_\beta} \e_{{\bf j}_\beta}+$ other terms. By equating the coefficients we take that $c_{1,0,0,0}=\psi_a$, hence
$c_{1,0,0,0}$ is uniquely determined. Now, the maximum term appeared into $R_1$ is $\de_{{\bf i}_\beta} \e_{{\bf j}_\beta}$. A second visit at the REM-Step 3e will give
$R_1=R_2+c_{2,0,0,0}\de_{{\bf i}_{2,0,0,0}}\e_{{\bf j}_{2,0,0,0}}$$\ast[L_{2,0,0},M_{2,0,0}]$. Working exactly as before we shall finally get that
$c_{2,0,0,0}=\psi_\beta-c_{1,0,0,0}(s_{1,0,0,j_{m,a}-1})$. 
Since
this expression is a function of $c_{1,0,0,0}$ and of some parameters
$w_{ijhk},s_{ijhk}$ considered as constants, we conclude that $c_{2,0,0,0}$ is
also defined uniquely. By induction, we can then see that all the coefficients are
uniquely determined. The polynomials $R_{\de,h},R_{\e,h},R_{\de\e,h^*}$ consist only of zero terms. These terms arise
either from the polynomial $A$ or from the products $c_{\lambda,u,\omega,h}\de_{{\bf i}_{\lambda,u,\omega,h}}\e_{{\bf j}_{\lambda,u,\omega,h}}$.
 The unique determination of the coefficients
 $c_{\lambda,u,\omega,h}$ entails the uniqueness of the remainders and, thus, the theorem has been proved.$\Box$
 \begin{rem}
 The equation (\ref{FLFMAIN}) can be written shortly as:
\[ \sum_{\bf \theta} \de_0^{\kappa_{\bf \theta}}\e_0^{\sigma_{\bf \theta}}\cdot (c_{\bf \theta}\de_{{\bf i}_{\bf \theta}}\e_{{\bf j}_{\theta}} \ast [L_{\bf \theta},M_{\bf \theta}])+{\bf R}_{\de,h}+{\bf R}_{\e,h}+R_{\de\e}\]
with ${\bf \theta}$ is an abbreviation for the vector: $(\lambda,u,\omega,h)$. 
\end{rem}

\noindent
Let $A$ be a $\de \e$-polynomial, ${\cal L,M}$, its sets of the Formal Linear Factors.
and ${\cal C}$ the set of the companion coefficients. The elements of these sets, or of any subset of them, can be transformed to specific $\de$ or $\e$-linear polynomials or to specific companion coefficients, by giving to the parameters
$w_{\lambda,u,h,k}$, $s_{\lambda,u,h,k}$ concrete values, according to a set of substitution rules. In this case we say that ${\cal L,M}$ and ${\cal C}$ are evaluated over a set of rules $U$, thus writing: $\left.
\begin{array} {c}
{\cal L}\\
\end{array} \right|_{U} $,
$\left.
\begin{array} {c}
{\cal M}\\
\end{array} \right|_{U} $ and
$\left.
\begin{array} {c}
{\cal C}\\
\end{array} \right|_{U} $.
More rigorously, let ${\bf W}=( w_{\lambda,i,h,k})$ be the set of the $w$-parameters, written as
a vector and ${\bf r}=( a_{\lambda,i,h,k}
)$ a vector of real numbers, which is in one-to-one
correspondence with the vector ${\bf W}$. We say that these parameters follow the
rule ${\bf r}$, thus writing ${\bf W \to r}$, if the following
substitutions are valid $ w_{\lambda,i,h,k}=a_{\lambda,i,h,k}$.
The set $\left.
\begin{array} {c}
{\cal L}\\
\end{array} \right|_{\bf r} $ is defined as follows:
\[ \left. \begin{array}
{c}
{\cal L}\\
\end{array} \right|_{\bf r} =\{ L_{\lambda,u,h}=\sum_{k=0}^{a_\lambda-1}w_{\lambda,u,h,k}\de_k+\de_{a_\lambda} \quad , \quad h=0,\ldots,h^*,
u=1,\ldots,u^*_h,\]
\[\lambda=1,\ldots,\lambda^*_{u,h} \mbox{ and } a_\lambda \mbox{ positive integers }
 , \quad \mbox{with} \quad {\bf W}=(w_{\lambda,u,h,k}) \to {\bf r} \} \]
Let $N$ be a set of
rules $N=\{ {\bf r}_1,{\bf r}_2, \ldots, {\bf r}_{\psi}\}$ then
$\left.
\begin{array} {c}
{\cal L}\\
\end{array} \right|_{N} $$=\{
\left.
\begin{array} {c}
{\cal L}\\
\end{array} \right|_{{\bf r}_1}, $
$\left.
\begin{array} {c}
{\cal L}\\
\end{array} \right|_{{\bf r}_2}, $
$\ldots ,$
$\left.
\begin{array} {c}
{\cal L}\\
\end{array} \right|_{{\bf r}_\psi}\} $
Simultaneously, if 
${\bf S}=( s_{\lambda,j,h,k})$ is the vector of the $s$-parameters and
${\bf q}=( b_{\lambda,j,h,k} )$ a vector of real numbers,
 which is in one-to-one
correspondence with the vector ${\bf S}$. We say that these parameters follow the
rule ${\bf q}$, thus writing ${\bf S \to q}$, if the following
substitutions are valid $ s_{\lambda,j,h,k}=b_{\lambda,j,h,k}$.
We define 
\[ \left. \begin{array}
{c}
{\cal M}\\
\end{array} \right|_{\bf q} =\{ M_{\lambda,\omega,h}=\sum_{k=0}^{b_\lambda-1}s_{\lambda,\omega,h,k}\e_k+\e_{b_\lambda} \quad , \quad h=0,\ldots,h^*,
\omega=1,\ldots,\omega^*_h,\]
\[\lambda=1,\ldots,\lambda^*_{\omega,h} \mbox{ and } b_\lambda \mbox{ positive integers }
 , \quad \mbox{with}, \quad {\bf S}=(s_{\lambda,\omega,h,k}) \to {\bf q} \} \]
Let $Q$ be the set of rules
$Q=\{ {\bf q}_1,{\bf
q}_2, \ldots, {\bf q}_{\mu}\}$, then
$\left.
\begin{array} {c}
{\cal M}\\
\end{array} \right|_{Q} $$=\{
\left.
\begin{array} {c}
{\cal M}\\
\end{array} \right|_{{\bf q}_1}, $
$\left.
\begin{array} {c}
{\cal M}\\
\end{array} \right|_{{\bf q}_2}, $
$\ldots ,$
$\left.
\begin{array} {c}
{\cal M}\\
\end{array} \right|_{{\bf q}_\mu}\} $. If $w$ and $s$-parameters appear together, then we say that they follow the rule $({\bf r,q})$, thus writing
$({\bf W,S})\to ({\bf r,q})$, if
$ w_{\lambda,i,h,k}=a_{\lambda,i,h,k}$ and 
$ s_{\lambda,j,h,k}=b_{\lambda,j,h,k}$. Let $({\bf r,q})$ be a concrete rule, then
$\left.
\begin{array} {c}
{\cal C}\\
\end{array} \right|_{\bf (r,q)}=$
$\{ c_{\lambda,u,\omega,h}(w_{\lambda,i,h,k},s_{\lambda,j,h,k})$
with
$({\bf W,S})\to ({\bf r,q})\}$.
Let $N,Q$ be two sets of rules  and $N\times Q$ their cartesian product, then
${\cal C}|_{N\times Q}=\{{\cal C}|_{({\bf r}_a,{\bf q}_b)}$ with $ ({\bf r}_a,{\bf q}_b) \in N\times Q\}$.  It is obvious that we can apply the above evaluation procedure not only to the entire sets ${\cal L,M}$ and ${\cal C}$ but to any subset of them, even to just one element. For instance, $L_{\lambda,u,h}|_{{\bf r}_1}$ denotes the evaluation of the 
$\de$-polynomial $L_{\lambda,u,h}$ over the specific substitution rule ${\bf r}_1$.
\begin{exa}\label{exaLFL}

Let us have the $\de\e$-polynomial:
\[ A=\de_0-2\de_1+2\de_1\de_2+4\de_2^2+\frac{1}{2}\de_2\de_3-2\de_2\de_4+\frac{1}{4}\de_4^2+3\e_2+16\e_1^2-18\e_1\e_2-2\e_1^2\e_2+\]
\[+5\e_2^2+2\e_1\e_2^2-\frac{1}{2}\e_2^3+\e_1^2\e_3-\e_1\e_2\e_3+\frac{1}{4}\e_2^2\e_3+\de_1\de_2\e_1+\frac{1}{4}\de_2^2\e_1+\frac{1}{2}\de_1\de_3\e_1-\]
\begin{equation}\label{para-de}
-\frac{1}{16}\de_3^2\e_1-4\de_1\de_2\e_2+\de_2^2\e_2+2\de_1\de_3\e_2-\frac{1}{4}\de_3^2\e_2
\end{equation}
In order to clarify the notion of Linear Factors, we shall follow the FLF-Algorithm step by step. We shall begin by working with the $\de$-part of $A$,
$A_{\de}=\de_0-2\de_1+2\de_1\de_2+4\de_2^2+$$\frac{1}{2}\de_2\de_3-2\de_2\de_4+\frac{1}{4}\de_4^2$. The linear part of $A_{\de}$ is $\de_0-2\de_1$ and hence $L_{0,1,0}=\de_0-2\de_1$. This polynomial, which is not formal, is the first member of the set ${\cal L}$. Now, the maximum term of the nonlinear part of $A_{\de}$ is $\frac{1}{4}\de_4^2$ and the REM-Substep 3e will give: 
$R=A_{\de,nl}-\frac{1}{4}\de_0^2\e_e \ast [L_{1,1,0},M_{1,0,0}]$ with
$L_{1,1,0}=w_{1,1,0,0}\de_0+$
$w_{1,1,0,1}\de_1+$
$w_{1,1,0,2}\de_2+$
$w_{1,1,0,3}\de_3+\de_4$, $M_{1,0,0}=\e_e$, this $L_{1,1,0}$ is the next member of ${\cal L}$ and $ccoef(L_{1,1,0})=-\frac{1}{4}$, ${\cal O}={\cal O} \cup \{\de_0^2 \e_e\}$. Now, the higher ordered term of the obtained polynomial $R$, is 
$-\frac{w_{1,1,0,3}}{2}\de_3\de_4$ and the REM-Substep 3e will give
$R_{new}=R_{old}+\frac{w_{1,1,0,3}}{2}\de_0\de_1\e_e \ast [L_{2,1,0}, M_{2,0,0}]$
with
$L_{2,1,0}=w_{2,1,0,0}\de_0+$
$w_{2,1,0,1}\de_1+$
$w_{2,1,0,2}\de_2+$
$\de_3$ and $M_{2,0,0}=\e_e$. Working similarly we shall finally get the next epression for $A_{\de}$:
\[ A_{\de}=\de_0-2\de_1+\frac{1}{4}\de_0^2 \ast (w_{1,1,0,0}\de_0+w_{1,1,0,1}\de_1+
w_{1,1,0,2}\de_2+w_{1,1,0,3}\de_3+\de_4)-\frac{w_{1,1,0,3}}{2}\de_0\de_1\ast\]
\[ \ast ( w_{2,1,0,0}\de_0+w_{2,1,0,1}\de_1+w_{2,1,0,2}\de_2+\de_3)-\left(-2-\frac{w_{1,1,0,2}}{2}+\frac{w_{1,1,0,3}w_{2,1,0,2}}{2}\right) \de_0\de_2 \ast\]  
\[ \ast (w_{3,1,0,0}\de_0 +w_{3,1,0,1}\de_1+\de_2)+c_4 \de_0\de_4 \ast (w_{4,1,0,0}\de_0+\de_1)+c_5\de_0^2\ast (w_{5,1,0,0}\de_0+w_{5,1,0,1}\de_1+w_{5,1,0,2}\de_2+\de_3)+\]
\[+c_6\de_0\de_1 \ast(w_{6,1,0,0}\de_0+w_{6,1,0,1}\de_1+\de_2)+c_7\de_0\de_2\ast (w_{7,1,0,0}\de_0+\de_1)+c_8\de_0^2 \ast (w_{8,1,0,0}\de_0+w_{8,1,0,1}\de_1+\de_2)+\]
\[+c_9 \de_0\de_1 \ast (w_{9,1,0,0}\de_0 +\de_1)+c_{10} \de_0^2 \ast (w_{10,1,0,0}\de_0+\de_1) +R\]
where $R$ contains only zero terms, particularly
$R=r_0\de_0^2+r_1\de_0\de_1+r_2\de_0\de_2 +r_3\de_0\de_3$. The coefficients $c_4,c_5,c_6,c_7,c_8,c_9,c_{10}$,$r_0,r_1,r_2,r_3$ are polynomial functions of the $w_{\lambda,u,h,i}$ parameters. They are not presented explicitly due to their large size (the coefficient $c_{10}$, for instance, contains 77 monomials). We factorize $R$ by $\de_0$ and we get:
$R=\de_0 \tilde{R}=\de_0 ( r_0\de_0+r_1\de_1+r_2\de_2 +r_3\de_3)$ and we repeat the procedure with
$A_{\de}=r_0\de_0+r_1\de_1+r_2\de_2 +r_3\de_3$. Since it is linear the SUBROUTINE REMAINDER will terminate at REM-STEP 2, giving $L_{0,2,0}=r_0\de_0+r_1\de_1+r_2\de_2 +r_3\de_3$ and
$ccoef(L_{0,2,0})=1$, ${\cal O}={\cal O} \cup \{\de_0\}$. All the above linear polynomials are elements of ${\cal L}$.
We are now in Step 4, working with the $\e$-part of $A$, $A_{\e}=3\e_2+16\e_1^2-$$18\e_1\e_2-2\e_1^2\e_2+$
$5\e_2^2+2\e_1\e_2^2-$$\frac{1}{2}\e_2^3+\e_1^2\e_3-$$\e_1\e_2\e_3+\frac{1}{4}\e_2^2\e_3$.
Initially, we get $M_{0,1,0}=3\e_2$ and for the nonlinear part of $A_{\e}$ we shall finally get:
\[ A_{\e,nl}=\frac{1}{4}\e_0^2\e_1\ast (s_{1,1,0,0}\e_0+s_{1,1,0,1}\e_1+\e_2)+\left(-1-\frac{s_{1,1,0,1}}{2}\right)\e_0\e_1\e_2\ast (s_{2,1,0,0}\e_0+\e_1)+\]
\[+\left(1-\frac{s_{1,1,0,1}}{4}+s_{2,1,0,0}+\frac{s_{1,1,0,1}s_{2,1,0,0}}{2}\right)\e_0^2\e_2\ast (s_{3,1,0,0}\e_0+\e_1)+\left(-\frac{1}{2}-\frac{s_{1,1,0,1}}{4}\right)\e_0^3 \ast\]
\[\ast (s_{4,1,0,0}\e_0+s_{4,1,0,1}\e_1+\e_2)+\tilde{c}_5\e_0\e_1^2\ast (s_{5,1,0,0}\e_0+\e_1)+\tilde{c}_6\e_0^2\e_1 \ast (s_{6,1,0,0}\e_0+\e_1)+\]
\[+\tilde{c}_7\e_0^3 \ast (s_{7,1,0,0}\e_0+\e_1)+5\e_0^2\ast (s_{8,1,0,0}\e_0+s_{8,1,0,1}\e_1+\e_2)+(-18-10s_{8,1,0,1})\e_0\e_1\ast\]
\[*(s_{9,1,0,0}\e_0+\e_1)+\tilde{c}_{10}\e_0^2 \ast (s_{10,1,0,0}\e_0+\e_1)+R\]
The coefficients $\tilde{c}_5,\tilde{c}_6,\tilde{c}_7,\tilde{c}_{10}$ have large size and therefore are not written in details. The quantity $R$ is
$R=r_{2,3}\e_0\e_2\e_3+$
$r_{1,2}\e_0\e_1\e_2+$
$r_{0,3}\e_0^2\e_3+$
$r_{0,1}\e_0^2\e_1+$
$r_{2,2}\e_0\e_2^2+$
$r_{1,3}\e_0\e_1\e_3+$
$r_{0,0}\e_0^3+r_{1,1}\e_0\e_1^2$
$+r_{0,2}\e_0^2\e_2$
$+r_2\e_0\e_2+r_1\e_0\e_1+r_0\e_0^2$. Factorizing $R$ by $\e_0$ we get $R=\e_0 \cdot \tilde{R}$ and then we repeat the substeps of step 4 once more, taking $\omega=2$, $M_{0,2,0}=r_{0}\e_0+r_{1}\e_1+r_{2}\e_2$ and 
$R=r_{2,3}\e_2\e_3+$
$r_{1,2}\e_1\e_2+$
$r_{0,3}\e_0\e_3+$
$r_{0,1}\e_0\e_1+$
$r_{2,2}\e_2^2+$
$r_{1,3}\e_1\e_3+$
$r_{0,0}\e_0^2+r_{1,1}\e_1^2$
$+r_{0,2}\e_0\e_2$. The running of the substeps of REM-STEP 3 will finally give
$R=\tilde{r}_{0,1}\e_0\e_1 \ast (s_{1,2,0,0}\e_0+s_{1,2,0,1}\e_1+\e_2)+$
$\tilde{r}_{0,2}\e_0\e_2 \ast (s_{2,2,0,0}\e_0+\e_1)$
$+\tilde{r}_{0,0}\e_0^2 \ast (s_{3,2,0,0}\e_0 + s_{3,2,0,1}\e_1+\e_2)$
$+\tilde{r}_{0,1}\e_0\e_1 \ast (s_{4,2,0,0}\e_0+\e_1)$
$+\tilde{r}_{0,0}\e_0^2\ast (s_{5,2,0,0}\e_0+\e_1)+R$
with $R=\hat{r}_0\e_0^2 +\hat{r}_1\e_0\e_1+\hat{r}_2\e_0\e_2+\hat{r}_3\e_0\e_3$ 
and $\tilde{r}_{0,0},\tilde{r}_{0,1},\tilde{r}_{0,2}$,
$\hat{r}_0,\hat{r}_1,\hat{r}_2,\hat{r}_3$ are polynomials of all the above parameters. Visiting SUB-Step 4b once more we shall take $R=\e_0\cdot (\hat{r}_0\e_0 +\hat{r}_1\e_1+\hat{r}_2\e_2+\hat{r}_3\e_3)$. Now, $\omega=3$ and since $\tilde{R}$ is linear
the SUBROUTINE REMAINDER will terminate by giving 
$M_{0,3,0}=\hat{r}_0\e_0 +\hat{r}_1\e_1+\hat{r}_2\e_2+\hat{r}_3\e_3$. All the above appeared linear polynomials are elements of ${\cal M}$.
We go now to step 5, working with the $\de\e$-part of $A$, $A_{\de\e}=\de_1\de_2\e_1+\frac{1}{4}\de_2^2\e_1+$
$\frac{1}{2}\de_1\de_3\e_1$
$-\frac{1}{16}\de_3^2\e_1-4\de_1\de_2\e_2+\de_2^2\e_2+$
$2\de_1\de_3\e_2-\frac{1}{4}\de_3^2\e_2$. We have $u=0,\omega=0$ and successive applications of the substeps of REM-STEP 3, will finally give the following expression for $A_{\de\e}$:
\[ A_{\de\e}=-\frac{1}{4}\de_0^2\e_0\ast [w_{1,0,0,0}\de_0+w_{1,0,0,1}\de_1+w_{1,0,0,2}\de_2+\de_3,s_{1,0,0,0}\e_0+s_{1,0,0,1}\e_1+\e_2]+\]
\[+\left(-\frac{1}{16}+\frac{s_{1,0,0,1}}{4}\right)\de_0^2\e_0 \ast [w_{2,0,0,0}\de_0+w_{2,0,0,1}\de_1+w_{2,0,0,2}\de_2+\de_3,s_{2,0,0,0}\e_0+\e_1]+\]
\[+\left(\frac{s_{1,0,0,0}}{4}+\frac{s_{2,0,0,0}}{16}-\frac{s_{1,0,0,1}s_{2,0,0,0}}{4}\right)\de_0^2\e_0 \ast [w_{3,0,0,0}\de_0+w_{3,0,0,1}\de_1+w_{3,0,0,2}\de_2+\de_3,\e_0]+\]
\[+\frac{w_{1,0,0,2}}{2}\de_0\de_1\e_0 \ast [w_{4,0,0,0}\de_0+w_{4,0,0,1}\de_1+\de_2,s_{4,0,0,0}\e_0+s_{4,0,0,1}\e_1+\e_2]+p_5\de_0\de_1\e_0 \ast[w_{5,0,0,0}\de_0+w_{5,0,0,1}\de_1+\]
\[+\de_2,s_{5,0,0,0}\e_0+\e_1]+p_6\de_0\de_1\e_0\ast [w_{6,0,0,0}\de_0+w_{6,0,0,1}\de_1+\de_2,\e_0]+\]
\[+\left(2+\frac{s_{1,0,0,1}}{2}-\frac{s_{1,0,0,2}s_{4,0,0,1}}{2}\right)\de_0\de_2\e_0 \ast [w_{7,0,0,0}\de_0+\de_1, s_{7,0,0,0}\e_0+s_{7,0,0,1}\e_1+\e_2]+\]
\[+p_8\de_0\de_2\e_0 \ast [w_{8,0,0,0}\de_0+\de_1, s_{8,0,0,0}\e_0+\e_1]+p_9\de_0\de_2\e_0 \ast [w_{9,0,0,0}\de_0+\de_1, \e_0]+p_{10}\de_0^2\e_0 \ast \]
\[ \ast [w_{10,0,0,0}\de_0+w_{10,0,0,1}\de_1+\de_2,s_{10,0,0,0}\e_0+s_{10,0,0,1}\e_1+\e_2]+p_{11}\de_0^2\e_0 \ast [w_{11,0,0,0}\de_0+w_{11,0,0,1}\de_1+\de_2,\]
\[ s_{11,0,0,0}\e_0+\e_1]+p_{12}\de_0^2\e_0 \ast [w_{12,0,0,0}\de_0+w_{12,0,0,1}\de_1+\de_2,\e_0]+p_{13}\de_0\de_1\e_0\ast[w_{13,0,0,0}\de_0+\de_1,\e_2]+\]
\[+p_{14}\de_0\de_1\e_0\ast[w_{14,0,0,0}\de_0+\de_1,s_{14,0,0,0}\e_0+\e_1]+\cdots +R\]
where $R=q_0\de_0^2\e_0+q_1\de_0\de_1\e_0+q_2\de_0\de_2\e_0+q_3\de_0\de_3\e_0$. Since $R=\de_0\e_0 \cdot \tilde{R}$
$=\de_0\e_0 (q_0\de_0+q_1\de_1+q_2\de_2+q_3\de_3)$ the algorithm will terminate with $L_{0,0,1}=q_0\de_0+q_1\de_1+q_2\de_2+q_3\de_3$. Summarizing, the final output of the algorithm is
${\cal L}=\{\de_0-2\de_1, w_{1,1,0,0}\de_0+w_{1,1,0,1}\de_1+w_{1,1,0,2}\de_2+w_{1,1,0,3}\de_3+\de_4, w_{2,1,0,0}\de_0+w_{2,1,0,1}\de_1+w_{2,1,0,2}\de_2+\de_3, \ldots \}$, ${\cal M}=\{3\e_2,s_{1,1,0,0}\e_0+s_{1,1,0,1}\e_1+\e_2,s_{2,1,0,0}\e_0+\e_1,\ldots\}$,
${\cal O}=\{\de_0,\de_4^2,\de_0\de_1,\ldots,\e_0,\e_0^2\e_1,\ldots,\de_0^2\e_0,\ldots\}$,
${\cal C}=\{1,\frac{1}{4},-\frac{w_{1,1,0,3}}{2}, -2-\frac{w_{1,1,0,2}}{2}+\frac{w_{1,1,0,3}w_{2,1,0,2}}{2}$,
$c_4,c_5,\ldots, \frac{1}{4},-1-\frac{s_{1,1,0,1}}{2},\ldots,\tilde{c}_5,\ldots$
$-\frac{1}{4},-\frac{1}{16}+\frac{s_{1,0,0,1}}{4},\ldots,p_5,p_5,\ldots\}$ and $R_{\e,\nu}=R_{\de,\nu}=0$, $h^*=1, R_{\de\e,0}=R_{\de\e,1}=0$. Moreover,
${\cal L}_{\de}={\cal L}-\{r_0\de_0+r_1\de_1+r_2\de_2+r_3\de_3\}$,
${\cal L}^*={\cal L}-\{\de_0-2\de_1\}$, ${\cal M}^*={\cal M}-\{3 \e_2\}$, 
$\overline{\cal L}={\cal L}-\{w_{1,0,0,0}\de_0+w_{1,0,0,1}\de_1+w_{1,0,0,2}\de_2+\de_3, \ldots\}$,
$\overline{\cal M}={\cal M}-\{s_{1,0,0,0}\e_0+s_{1,0,0,1}\e_1+\e_2, \ldots \}$,
$\overline{ \cal L}^*=\overline{ \cal L}-\{\de_0-2\de_1\}$, $\overline{\cal M}^*=\overline{\cal M}-\{3\e_2 \}$. Let us consider now the rules ${\bf r}_1=\{ w_{h,u,\omega,\lambda}=1\}$, ${\bf r}_2=\{ w_{1,1,0,3}=0,w_{1,1,0,2}=-4,w_{3,1,0,1}=0,$
$w_{4,1,0,0}=0,w_{1,1,0,1}=0,$ all the other w-parameters take arbitray values $\}$, ${\bf q}_1=\{s_{h,u,\omega,\lambda}=1\}$, ${\bf q}_2=\{s_{1,0,0,1}=0 $ all the other s-parameters take arbitray values $\}$, then ${\cal L}|_{{\bf r}_1}=\{\de_0-2\de_1,\de_0+\de_1+\de_2+\de_3+\de_4,\de_0+\de_1+\de_2+\de_3,\ldots\}$,
${\cal L}_{\de}|_{{\bf r}_2}=\{\de_0-2\de_1,w_{1,1,0,0}\de_0-4\de_2,w_{2,1,0,0}\de_0+w_{2,1,0,1}\de_1+w_{2,1,0,2}\de_2+\de_3, \ldots \}$, $A_{\de}|_{{\bf r}_2}=\de_0-2\de_1+\frac{1}{4}\de_0^2 \ast (w_{1,1,0,0}\de_0-4\de_2)+\frac{1}{2}\de_0\de_1 \ast (w_{6,1,0,0}\de_0+w_{6,1,0,0}\de_1+\de_2)+R|_{{\bf r}_2}$, ${\cal M}|_{{\bf q}_1}=\{3\e_2,\e_0+\e_1+\e_2,\e_0+\e_1,\ldots \}$ and finally
$c_{2,0,0,0}|_{({\bf r}_2,{\bf q}_2)}=\left(-\frac{1}{16}+\frac{s_{1,0,0,1}}{4}\right)|_{({\bf r}_2,{\bf q}_2)}=-\frac{1}{16}$.

\end{exa}
\section{The Model Matching Problem}
In this section we present the main application of our methodology.
Before we examine it, 
we would like to begin by a short algebraic
description of the non-linear discrete input-output systems, 
via the notion of $\de
\e$-polynomials.
Suppose that we have a system of the form (\ref{main}). By using
$\de\e$-operators we can rewrite this as follows: $F[y(t),u(t)]=0$, where $F$ is a $\de\e$-polynomial. As usual, many times we decompose $F$ to its pure $\de$, $\e$ and $\de\e$-parts, i.e. $F=F_{\de}+F_{\e}+F_{\de\e}$ and $F_{\de},F_{\e}$ to their linear and nonlinear parts: $F_{\de}=F_{\de,l}+F_{\de,nl}$, $F_{\e}=F_{\e,l}+F_{\e,nl}$. ($F_{\de\e}$ is already purely nonlinear). The causality and
solvability of the system is guaranteed by the equality $d(F)=d(F_{\de,l})$.
Often, we shift all the delays properly so that $d(F)=d(F_{\de,l})=0$,
 which means  
that the lower delayed term, i.e. $y(t)$, appears
in the linear $\de$-part of the system and thus we can solve (\ref{main}) with respect
to $y(t)$ in a direct way. To each nonlinear system of the form (\ref{main}) we assign a vector
of real numbers ${\bf y}_0=(y_0,y_1, \ldots, y_{k-1})$ which
gives the so called initial conditions:
$y(0)=y_0$, $y(1)=y_1$, \ldots, $y(k-1)=y_{k-1}$, where $k$ is the maximum
delay appeared in the output signal. Since the signals, involved in
(\ref{main}), are causal, i.e. $y(t)=0$, $u(t)=0$, for $t<0$, for each
given vector of initial conditions any input signal $u(t)$ determines
 a unique output signal that satisfies (\ref{main}), for $n \ge k$. Let $K=\sum_{i=0}^{k-1}\tau_i\de_i$ be a linear $\de$-polynomial, we say that the initial conditions ${\bf y}_0$ satisfy the relation $K{\bf y}_0=0$ if $\sum_{i=0}^{k-1}\tau_iy_{k-1-i}=0$. 
Let us have two causal nonlinear systems: $F[y(t),u(t)]$
with ${\bf y}_0=(y_0,y_1, \ldots, y_{k-1})$ and
$\hat{F}[\psi(t),v(t)]$ with 
${\bf \psi}_0=(\psi_0,\psi_1,\ldots,\psi_{\lambda-1})$ and $k>\lambda$.
We say that the two systems " operate "
under identical initial conditions,
if $\psi_0=y_0$,
$\psi_1=y_1$, $\ldots$, $\psi_{\lambda-1}=y_{\lambda-1}$ 
and $y_{\lambda}=\psi(\lambda)$, 
$y_{\lambda+1}=\psi(\lambda+1)$, $\ldots$, 
$y_{k-1}=\psi(k-1)$. In other words we must give
 as initial conditions to the system which starts to product outputs
 later, the corresponding outputs of the other system.
The next theorem is crucial for the model mathcing problem. Its proof can be found at \cite{kn:difeq}.
\begin{theo}\label{eksdiaf} 
Let $Gy(t)=0$ and $\Phi y^*(t)=0$, $G,\Phi$ a nonlinear and a linear $\de$-polynomials\, be two systems without inputs, a nonlinear and a linear one, (actually they are difference equations with respect to $y(t)$ and $y*(t)$ ). If $G=\hat{G}\ast \Phi$, then $y(t)=y^*(t)$ whenever the two systems operate under identical initial conditions.
If $G=\hat{G}\ast \tilde{\Phi}$ and $\Phi=\hat{\Phi} \ast \tilde{\Phi}$, then $y(t)=y^*(t)$ whenever the two systems operate under identical initial conditions ${\bf y}_0$, which satisfy the relation $\tilde{\Phi}{\bf y}_0=0$.
\end{theo}
\noindent
The model matching problem for input - output discrete systems has already attracted a
lot of attention. Though there is already a clear understanding of the linear problem
\cite{kn:gladbook}, \cite{kn:amstrong}, 
 the non-linear case is still being examined
and is the subject matter of a large number of papers (see for instance
\cite{kn:phd1},\cite{kn:yamanaka},\cite{kn:gladmodelreference}, to mention but a few).
\par
\noindent In this paper we examine the Model Matching Problem (MMP) for non-linear
discrete systems of the form (\ref{main}), or $F[y(t),u(t)]=0$, $F$ a $\de \e$-polynomial, by using the entire algebraic background
that has been previously developed. This approach is a continuation  of research on
non-linear discrete systems in \cite{kn:bibo},\cite{kn:phd1},\cite{kn:factn},\cite{kn:kotritt},
\cite{kn:ecc05}. Our aim is to find a feedback
relation of the form $u(t)=Sy(t)$, $S$ a linear $\de$-polynomial, so that if it is fed
back to the original system we get an output that is identical
 to that of a given
desired linear system, described by the equation $A_dy^*(t)=0$, $A_d$ a linear
$\de$-polynomial. If this holds, we say that the nonlinear systems matches the linear one.
\par
\noindent
To achieve a solution of the problem we work as follows: First by means 
of the $FLF$-Algorithm we obtain all the linear factors of the nonlinear
system and then we 
rewrite $F$ as in (\ref{FLFMAIN}). By substituting the feedback relation $u(t)=Sy(t)$, we eliminate the input signal $u(t)$ and we take an expression containing only the output signal $y(t)$, actually we have a difference equation with respect to $y(t)$. The majority of the terms, appeared in this difference expression contain either the formal $\de$-linear factors or the formal $\e$-linear factors of $F$ or the " free " linear polynomial $S$. The cornerstone of our approach is the following: if we can achieve to find values of the parameters $w_{ijhk},s_{ijhk}$ and of the free linear polynomial $S$, so that all the terms of the said difference equation (formal and no-formal), have a common factor which is a factor of $A_d$ as well, then, these specific values of $S$ will constitute a solution of the model matching problem. To accomplish this task we examine several cases. In all of them we eliminate the unneccesary remainders and some of the linear factors and then we " ask " from the linear parts of the original systems (which do not have parametrical coefficients and thus cannot be influenced by the choice of proper values) to have a common factor. This can be achieved by solving certain Diophantine equations.  Two things are remarkable here, first we do not get a single solution but a class of solutions, parameterized
via a linear polynomial $Q$. Proper choice of this polynomial
can create a causal feedback-law or a law with as much or as few delays
we like. Secondly, and most important, all the calculations and procedures are devoted to the star-product among linear polynomials: the linear factors of $F$. But, as we have seen in section \ref{prwto} the star-product of linear $\de$ or $\e$-polynomials coincides with the usual product among univariate polynomials. Hence, we can use all the available results and tools of the classical algebraic geometry or combutational algebra, to carry out the manipulations we need.
Before we present a symbolic algorithm, named the MM-Algorithm, which solves the problem completely, we shal describe the GCD-Values Algorithm. This algorithm will be used as a subroutine of the main algorithm. What it essentially does is to discover those values of the parameters which permit to a set of formal linear polynomials to have a common factor. To achieve that it works in two steps: firstly it discovers the values of the parameters which make certain companion coefficients no-zero and thus the corresponding linear factors are present in the expression (\ref{big}) and then it selects among these values, those specific ones which allow the corresponding linear polynomials to have a common factor.  
\small
{\sf \vskip 15 pt \noindent \underline{\bf The GCD-Values-Algorithm} 

\vskip 7 pt \noindent
{\bf Input:} The set a formal linear polynomials ${\cal F}$. 
\vskip 5 pt 
\noindent 
{\bf Output:} The set of rules $GCDV({\cal F})$.
\vskip 5 pt \noindent
{\bf STEP 0:} We set $GCDV({\cal F})=\{\}$ and $m=card({\cal F})$, that is $m$ is the amount of polynomials contained into ${\cal F}$.
\vskip 5 pt \noindent
{\bf STEP 1:} {\bf FOR} $n=1$ {\bf TO } $m$
\begin{quote}
{\bf Step 1a:} We construct the power set $\textbf{P}$, consisting from all the subsets of ${\cal F}$, with $n$ elements, i.e. $\textbf{P}=\{ {\cal A}: {\cal A}\subset {\cal F}$ with $card({\cal A})=n\}$.
 \vskip 5 pt
\noindent {\bf STEP 2a:} {\bf REPEAT} for all ${\cal A} \in \textbf{P}$.
\begin{quote}
 \vskip 5 pt
\noindent {\bf STEP 2a-I:} Construct the set ${\cal C}^{\cal A}$ of the companion coefficients of the elements of ${\cal A}$.
 \vskip 5 pt
\noindent {\bf STEP 2a-II:} We form the set of rules $NZ=\{{\bf r}: {\bf r}$ is defined by the substitution $(w_{ijhk},s_{ijhk})=
(a_{ijhk},b_{ijhk})$, $a_{ijhk}, b_{ijhk} \in {\bf R}$, such that $ccoef(J)\ne 0$, for $J \in {\cal A}$ and $ccoef(\Psi)=0$ for $\Psi \in {\cal F}-{\cal A}\}$.
\vskip 5 pt
\noindent {\bf STEP 2a-IIÉ:} We form the set of rules $CF \subset NZ$ such that $\gcd(J_i)\ne \de_0$,$J_i \in {\cal A}|_{NZ}$ and we set $GCDV({\cal F})=GCDV({\cal F})\cup CF$.
\end{quote}
\vskip 5 pt
\noindent {\bf END OF THE REPEAT}
\end{quote}
\vskip 5 pt
\noindent {\bf END OF THE FOR}
}
\normalsize
\begin{rem} 1) The calculation of the values of the parameters which ensure that $\gcd(J_i)\ne \de_0$, can be carried out by means of resultans theory \cite{kn:cox}.
\par
\noindent
2) Obviously, if $\gcd(J_i)=\de_0$, for $J_i$ belonging to a specific set of polynomials ${\cal A}$, then for any set ${\cal B}$ with ${\cal A} \subset {\cal B}$, we shall have $\gcd(J_i)=\de_0$, $J_i \in {\cal B}$, too. Therefore, we can avoid to check any set which includes ${\cal A}$. This comment can reduce the whole amount of calculations, significantly.
\end{rem}
We are now ready to state the main algorithm of the paper.

\small
{\sf \vskip 15 pt \noindent \underline{\bf The MM-Algorithm} 

\vskip 7 pt \noindent
{\bf Input:} The $\de\e$-polynomial $F$. 
\vskip 5 pt 
\noindent 
{\bf Output:} The sets ${\bf F_1,F_2,F_3}$.
\vskip 5 pt \noindent
\begin{quote}
{\bf STEP 1:} We decompose the polynomial $F$ into its linear and
non-linear parts: $F=F_{\de,l}+F_{\e,l}+F_{nl}$ and we rename $A=F_{\de,l}$ and $B=F_{\e,l}$.
 \vskip 5 pt
\noindent {\bf STEP 2:} By means of the $FLF$-Algorithm we form the 
sets ${\cal L}$, ${\cal M}$, ${\cal L}^*$, ${\cal M}^*$, ${\cal L}_{\de}$, 
 $\overline{\cal L}$, $\overline{\cal M}$, 
$\overline{\cal L}^*$,
$\overline{\cal M}^*$ and $R_{\de,h}=R_{\e,h}$, $h=0,\ldots,h^*$, $R_{\de\e}$.
\vskip 5 pt \noindent {\bf STEP 3:} By means of the GCD-Values-Algorithm, whenever it is needed, we form the following sets of rules:
\begin{itemize}
\item $LV=GCDV({\cal L}_{\de})$
\item $GV=GCDV({\cal L } \cup \overline{\cal M}) \cup GCDV(\overline{\cal L } \cup {\cal M})$
\item $MV=GCDV({\cal L }^* \cup \overline{\cal M}^*) \cup GCDV(\overline{\cal L }^* \cup {\cal M}^*)$
\item $DR=\{ (w_{ijhk})=
(a_{ijhk})$, $a_{ijhk} \in {\bf R}$, such that 
$R_{\de,0}=0\}$.
\item $RV=\{ (w_{ijhk},s_{ijhk})=
(a_{ijhk},b_{ijhk})$, $a_{ijhk}, b_{ijhk} \in {\bf R}$, such that 
$R_{\de,h}=R_{\e,h}=0$, $h=0,\ldots,h^*$, $R_{\de\e}=0 \}$

\end{itemize}
\vskip 5 pt \noindent {\bf STEP 4:} {\bf IF}
$ LV \cap DR \ne
\emptyset$ {\bf THEN}
put ${\bf F}_1=
$ $ \{\Phi : \Phi $ a common factor of $L_i \in {\cal L}_{\de}|_{LV\cap DR}\}$.

\vskip 5 pt \noindent {\bf STEP 5:}  {\bf IF}
$RV \cap GV \ne \emptyset$ {\bf THEN}
put ${\bf F}_2=
$ $ \{\Phi : \Phi $ a common factor of 
$ L_i \in
\left. \begin{array} {c}
{\cal L}\\
\end{array} \right|_{RV \cap GV}$ and $M_i \in \left. \begin{array}
{c}
\overline{\cal M}\\
\end{array} \right|_{RV \cap GV} \}$
$\cup$ 
$\{\Phi : \Phi $ a common factor of $  L_i \in
\left. \begin{array} {c}
\overline{\cal L}\\
\end{array} \right|_{RV \cap GV}$ and $M_i \in \left. \begin{array}
{c}
{\cal M}\\
\end{array} \right|_{RV \cap GV} \}$.

\vskip 5 pt \noindent {\bf STEP 6:}  {\bf IF}
$RV \cap MV \ne \emptyset$ {\bf THEN}
put ${\bf F}_3=$
$\{ \Phi: \Phi$ a common factor of $ L_i \in
\left. \begin{array} {c}
{\cal L}^*\\
\end{array} \right|_{RV \cap MV}$ and $M_i \in \left. \begin{array}
{c}
\overline{\cal M}^*\\
\end{array} \right|_{RV \cap MV} \}$
$\cup$
$\{ \Phi: \Phi=$ a common factor of $ L_i \in
\left. \begin{array} {c}
\overline{\cal L}^*\\
\end{array} \right|_{RV \cap MV}$ and $M_i \in \left. \begin{array}
{c}
{\cal M}^*\\
\end{array} \right|_{RV \cap MV}\}$

\vskip 5 pt \noindent {\bf STEP 7:}
Goto to the {\bf Output} 
\end{quote}
}
\normalsize
The next theorem includes the feedback construction upon request.

\begin{theo}\label{mainresult}
Let $F[y(t),u(t)]=0$, $F$ a $\de\e$-polynomial, be a nonlinear discrete input-output system with cross-products. Let $A$ be the linear output part of $F$, $-B$ its linear input (i.e. $F_{\de,l}=A,F_{\e,l}=-B$) and $A_dy^*(t)=0$, $A_d$ a linear $\de$-polynomial, a linear desired system. Let ${\bf F}_i, i=1,2,3$ be the outputs of the $MM$-Algorithm, after its application to the polynomial $F$. We suppose that at least one of the sets ${\bf F}_i, i=1,2,3$ is novoid. Then:
\par
a) If $A_d$ belongs to some of the sets ${\bf F}_i, i=1,2,3$, then the feedback law $u(t)=Sy(t)$, with
\[ S=
\left\{
\begin{array}{ll}
Q\ast A_d, & \mbox{if } A_d \in {\bf F}_1 \\
Q, & \mbox{if } A_d \in {\bf F}_2 \\
Z_0+Q\ast A_d, & \mbox{if } A_d \in {\bf F}_3 \\
\end{array}
\right.
\]
where $(R_0,Z_0)$ is a solution of the Diophantine equation $R\ast A_d+Z\ast B=A$, and $Q$ an arbitrary linear $\de$-polynomial, gives a closed-loop system which produces an output $y(t)$ equal to $y^*(t)$, whenever the nominal nonlinear system and the linear desired one, "operate" under identical initial conditions.
\par
b) If $A_d$ does not belong to any of the sets ${\bf F}_i, i=1,2,3$ but there are subsets
$\tilde{\bf F}_i \subseteq {\bf F}_i, i=1,2,3$, not all of them void, such that
$gcd(A_d,\Phi)=\tilde{\Phi}\ne \de_0$,for
$\Phi \in \tilde{\bf F}_i, i=1,2,3$, 
then, the feedback law $u(t)=Sy(t)$, with
\[ S=
\left\{
\begin{array}{ll}
Q\ast \tilde{\Phi}, & \mbox{if } \Phi \in \tilde{\bf F}_1 \\
Q, & \mbox{if } \Phi \in \tilde{\bf F}_2 \\
Z_0+Q\ast \tilde{\Phi}, & \mbox{if } \Phi \in \tilde{\bf F}_3 \\
\end{array}
\right.
\]
where $(R_0,Z_0)$ is a solution of the Diophantine equation $R\ast \tilde{\Phi}+Z\ast B=A$, and $Q$ an arbitrary linear $\de$-polynomial, gives a closed-loop system which produces an output $y(t)$ equal to $y^*(t)$, whenever the nominal nonlinear system and the linear desired one, "operate" under identical initial conditions ${\bf y}_0$, which satisfy the relation $\tilde{\Phi}{\bf y}_0=0$.

 \end{theo}
 
 {\bf Proof:} Let $F[y(t),u(t)]=0$ be the nominal plant, $F$ a $\de\e$-polynomial. We shall follow the MM-Algorithm step by step. By means of the FLF-Algorithm and theorem \ref{typos} we get the following expression for $F$.
\[  
 F=\sum_{h=0}^{h^*}\sum_{\omega=0}^{\omega_h^*}\sum_{u=0}^{u_h^*}\sum_{\lambda=0}^{v_{u,\omega,h}} \de_0^{\kappa_{u,h}}\e_0^{\sigma_{\omega,h}} \cdot (c_{\lambda,u,\omega,h}\de_{{\bf i}_{\lambda,u,\omega,h}}\e_{{\bf j}_{\lambda,u,\omega,h}}\ast [L_{\lambda,u,h},M_{\lambda,\omega,h}])+\]
 \begin{equation}\label{mm1}
 +{\bf R}_{\de,h}+{\bf R}_{\e,h}+R_{\de\e}
 \end{equation}
 Using the border values, given at \ref{typos}, and the fact that $L_{0,1,0}=F_{\de,l}=A$,
 $M_{0,1,0}=F_{\e,l}=-B$, we can split (\ref{mm1}) as follows:
 \[ F=A+\sum_{u=1}^{u_0^*}\sum_{\lambda=1}^{v_{u,0,0}}\de_0^{\kappa_{u,0}} \cdot c_{\lambda,u,0,0} \de_{{\bf i}_{\lambda,u,0,0}} \ast L_{\lambda,u,0}+{\bf R}_{\de,0}- B+\sum_{\omega=1}^{\omega_0^*}\sum_{\lambda=0}^{v_{0,\omega,0}}\e_0^{\sigma_{\omega,0}} \cdot c_{\lambda,0,\omega,0} \e_{{\bf j}_{\lambda,0,\omega,0}} \ast M_{\lambda,\omega,0}+\]

 \[+{\bf R}_{\e,0}+\sum_{\lambda=1}^{v_{0,0,0}}c_{\lambda,0,0,0}\de_{{\bf i}_{\lambda,0,0,0}}\e_{{\bf j}_{\lambda,0,0,0}} \ast [L_{\lambda,0,0},M_{\lambda,0,0}]+\de_0^{\kappa_{u_0^*,0}}\e_0^{\sigma_{\omega_0^*,0}} \sum_{h=1}^{h^*}\left(L_{0,1,h}+\sum_{u=1}^{u_h^*}\sum_{\lambda=1}^{v_{u,0,h}}\de_0^{\kappa_{u,h}}\right.\]
 
 \[ \cdot c_{\lambda,u,0,h} \de_{{\bf i}_{\lambda,u,0,h}} \ast L_{\lambda,u,h}+{\bf R}_{\de,h}+ +M_{0,1,h}+\sum_{\omega=1}^{\omega_h^*}\sum_{\lambda=1}^{v_{0,\omega,h}}\e_0^{\sigma_{\omega,h}} \cdot c_{\lambda,0,\omega,h} \e_{{\bf j}_{\lambda,0,\omega,h}} \ast M_{\lambda,\omega,h}+{\bf R}_{\e,h}+\]

\begin{equation}\label{mm2a}
\left.+\sum_{\lambda=1}^{v_{0,0,h}} \de_0^{\kappa_{u_h^*,h}}\e_0^{\sigma_{\omega_h^*,h}} \cdot c_{\lambda,0,0,h}\de_{{\bf i}_{\lambda,0,0,h}}\e_{{\bf j}_{\lambda,0,0,h}} \ast [L_{\lambda,0,h},M_{\lambda,0,h}]\right)+R_{\de\e}
\end{equation} 
 The above expression (\ref{mm2a}), recovers the " output " part, the " input " part and the " cross-products " part, contained into the original polynomial $F$. By substituting (\ref{mm2a}) and the feedback relation $u(t)=Sy(t)$ into the original nonlinear discrete system, we turn it to the following difference equation with respect to $y(t)$:
 \[ Ay(t)+\sum_{u=1}^{u_0^*}\sum_{\lambda=1}^{v_{u,0,0}}\de_0^{\kappa_{u,0}}y(t) \cdot c_{\lambda,u,0,0} \de_{{\bf i}_{\lambda,u,0,0}} \ast L_{\lambda,u,0}y(t)+{\bf R}_{\de,0}y(t)-B\ast Sy(t)+\sum_{\omega=1}^{\omega_0^*}\sum_{\lambda=1}^{v_{0,\omega,0}}\e_0^{\sigma_{\omega,0}}*Sy(t)\]

 \[ \cdot c_{\lambda,0,\omega,0} \e_{{\bf j}_{\lambda,0,\omega,0}} \ast M_{\lambda,\omega,0}\ast S y(t)+{\bf R}_{\e,0}\ast S y(t)+\sum_{\lambda=1}^{v_{0,0,0}}c_{\lambda,0,0,0}\de_{{\bf i}_{\lambda,0,0,0}}\e_{{\bf j}_{\lambda,0,0,0}} \ast [L_{\lambda,0,0}y(t),M_{\lambda,0,0}\ast Sy(t)]+\]

\[+\de_0^{\kappa_{u_0^*,0}}y(t)\e_0^{\sigma_{\omega_0^*,0}}\ast S y(t) \sum_{h=1}^{h^*}\left(L_{0,1,h}y(t)+\sum_{u=1}^{u_h^*}\sum_{\lambda=1}^{v_{u,0,h}}\de_0^{\kappa_{u,h}}y(t) \cdot c_{\lambda,u,0,h} \de_{{\bf i}_{\lambda,u,0,h}} \ast L_{\lambda,u,h}y(t)+ \right.\]

\[+{\bf R}_{\de,h}y(t)+M_{0,1,h} \ast Sy(t) +\sum_{\omega=1}^{\omega_h^*}\sum_{\lambda=1}^{v_{0,\omega,h}}\e_0^{\sigma_{\omega,h}}\ast Sy(t) \cdot c_{\lambda,0,\omega,h} \e_{{\bf j}_{\lambda,0,\omega,h}} \ast M_{\lambda,\omega,h} 
\ast S y(t)+{\bf R}_{\e,h}\ast S y(t)+\]

\[+\left.+\sum_{\lambda=1}^{v_{0,0,h}} \de_0^{\kappa_{u_h^*,h}}y(t)\e_0^{\sigma_{\omega_h^*,h}}\ast S y(t) \cdot c_{\lambda,0,0,h}\de_{{\bf i}_{\lambda,0,0,h}}\e_{{\bf j}_{\lambda,0,0,h}} \ast [L_{\lambda,0,h}y(t),M_{\lambda,0,h}\ast Sy(t)]\right)+\]
\begin{equation}\label{mm2}
R_{\de\e}\ast S [y(t),u(t)]=0
\end{equation} 
Where we used the fact that $A,L_{\lambda,u,h},R_{\de,h}$ are polynomials of $\de$-operators only, and thus they act exclusively on the output, whilst $B,M_{\lambda,u,h},R_{\e,h}$ are polynomials of $\e$-operators only, and thus they act exclusively on the input. By means of proposition \ref{star} (c) we can easily prove that:
\par
\noindent
(i) $\de_{\bf i} \e_{\bf j} \ast [Ly,M \ast S
y]=(\de_{\bf i} \ast L) \cdot (\de_{\bf j} \ast M \ast S) y(t) $, for any linear
polynomials $L,M,S$ (the operator $\e_j$ has changed to $\de_j$, since we are dealing
only with $y$ sequences).
\par
\noindent
(ii) $\e_0^\beta \ast Sy(t)=S^\beta y(t)$
\par
\noindent
(iii) $R_{\e,h} \ast Sy(t)=\sum_{\beta_h}c_{\beta_h}\e_0^{\beta_h}\ast Sy(t)=\sum_{\beta_h}c_{\beta_h}S^{\beta_h}y(t)$.
\par
\noindent
(iv) $R_{\de\e}=\sum_{(\alpha,\beta)}c_{(\alpha,\beta)}\de_0^{\alpha}\e_0^{\beta}\ast[y(t),Sy(t)]=\sum_{(\alpha,\beta)}c_{(\alpha,\beta)}\de_0^{\alpha}S^{\beta}y(t)$ $=y(t)]=\sum_{(\alpha,\beta)}c_{(\alpha,\beta)}y^a(t)[Sy(t)]^b$.
\par
\noindent
Therefore, after some manipulations, (\ref{mm2}) becomes:
\[ \left[(A-B*S)+\sum_{u=1}^{u_0^*}\sum_{\lambda=1}^{v_{u,0,0}}\de_0^{\kappa_{u,0}}\cdot c_{\lambda,u,0,0} \de_{{\bf i}_{\lambda,u,0,0}} \ast L_{\lambda,u,0}\right.+{\bf R}_{\de,0}+ \sum_{\omega=1}^{\omega_0^*}\sum_{\lambda=1}^{v_{0,\omega,0}} c_{\lambda,0,\omega,0} S^{\sigma_{\omega,0}} \cdot(\de_{{\bf j}_{\lambda,0,\omega,0}} \ast \]

\[\ast M_{\lambda,\omega,0}\ast S) +\sum_{\beta_0}c_{\beta_0}S^{\beta_0}+\sum_{\lambda=1}^{v_{0,0,0}}c_{\lambda,0,0,0}(\de_{{\bf i}_{\lambda,0,0,0}} \ast L_{\lambda,0,0}) \cdot ( \de_{{\bf j}_{\lambda,0,0,0}} \ast M_{\lambda,0,0}\ast S)+\]

\[+\de_0^{\kappa_{u_0^*,0}}S^{\sigma_{\omega_0^*,0}}\cdot \sum_{h=1}^{h^*}\left(L_{0,1,h}+\sum_{u=1}^{u_h^*}\sum_{\lambda=1}^{v_{u,0,h}}\de_0^{\kappa_{u,h}} \cdot c_{\lambda,u,0,h} \de_{{\bf i}_{\lambda,u,0,h}} \ast L_{\lambda,u,h}+{\bf R}_{\de,h}+ \right.\]
 
\[+M_{0,1,h} \ast S +\sum_{\omega=1}^{\omega_h^*}\sum_{\lambda=1}^{v_{0,\omega,h}} c_{\lambda,0,\omega,h} S^{\sigma_{\omega,h}} \cdot(\de_{{\bf j}_{\lambda,0,\omega,h}} \ast M_{\lambda,\omega,h} 
\ast S )+\sum_{\beta_h}c_{\beta_h}S^{\beta_h}+\]

\begin{equation}\label{momatcheq}
\left.+\sum_{\lambda=1}^{v_{0,0,h}}  c_{\lambda,0,0,h}   \de_0^{\kappa_{u_h^*,h}}S^{\sigma_{\omega_h^*,h}} \cdot  (\de_{{\bf i}_{\lambda,0,0,h}}\ast 
L_{\lambda,0,h})\cdot (\de_{{\bf j}_{\lambda,0,0,h}} \ast M_{\lambda,0,h}\ast S)\right)+\left.\sum_{(\alpha,\beta)}c_{(\alpha,\beta)}\de_0^{\alpha}S^{\beta}\right]y(t)=0
\end{equation}
The left-hand side of (\ref{momatcheq}) is a polynomial of $\de$-operators only. We denote it by $F_f$. So, (\ref{momatcheq}) is nothing else than the closed-loop system we take from $F[y(t),u(t)]=0$ by using the feedback-law $u(t)=Sy(t)$, in other words $F[y(t),Sy(t)]=F_fy(t)=0$. 
We are now ready to state the main argument of this proof. Since the desired system
$A_dy(t)=0$ has no input, it is obvious that the closed-loop system will have the same
dynamic behaviour as the desired one, only if we can find values  for the parametrical coefficients and the unknown
quantity $S$, such that the difference equations (\ref{momatcheq}) and $A_dy(t)=0$
 have common solutions.
To achieve this, we first need to use the equation (\ref{momatcheq}). By means of theorem \ref{eksdiaf}, 
we can establish that if we can find a set of values of the parameters so that $\Phi$ is a linear common factor of all the evaluated terms of (\ref{momatcheq}), then any solution of the difference equation $\Phi y(t)=0$, is also a solution of the equation (\ref{momatcheq}). If we ensure that this $\Phi$ is a factor of the polynomial $A_d$, too, or it coincides with it, then we have the desired result. Therefore, we
have to consider the following three cases, each of which corresponds to one of the three steps
4,5,6.
\par
\noindent 
{\bf Case I.} 
Let us suppose that ${\bf F}_1$ is non void and that $\Phi \in {\bf F}_1$. Step 4 of the MM-Algorithm implies that $LV\cap DV \ne \emptyset$. This condition means that we can find at least a rule of substitutions for the parameters $(w_{\lambda,u,h,k},s_{\lambda,\omega,h,k})$, named ${\bf q}$, such that two targets are accomplished simultaneously: (i) The remainder $R_{\de,0}$ is eliminated. This is due to the fact that ${\bf q} \in DR$. (ii) There is a certain amount of polynomials $L_{\lambda,u,0}$ with $ccoef(L_{\lambda,u,0})|_{\bf q}\ne 0$ (thus the polynomials $L_{\lambda,u,0}|_{\bf q}$ are not anihilated from $F_f|_{\bf q}$ and hence we can calculate their great common divider ) and $\gcd (L_{\lambda,u,0}|_{\bf q} \in {\cal L}_{\de}|_{\bf q})=\Phi\ne \de_0$. This is due to the fact that ${\bf q} \in LV$ and the construction of the GCDV-Algorithm.  Now, we set $S =Q \ast \Phi$, $Q$ an arbitrary linear polynomial. Since $A \in {\cal L}_{\de}|_{\bf q}$ , ($A$ has not parametrical coefficients) 
$\Rightarrow $
$A =T \ast \Phi$, for some $T$ and hence $A -B \ast S=T \ast \Phi - B \ast Q
\ast \Phi= (T-B \ast Q) \ast \Phi=K \ast \Phi$ and so $\Phi$ is a factor of $A-B\ast S$, too. Observing that  $S$ is a factor of all the terms of (\ref{momatcheq}), except the terms of the first line and taking into account all the above posed arguments and the structure of $F_f$, we finally conclude that $\Phi$ is a factor of all the terms of $F_f|_{\bf q}$, in other words $F_f|_{\bf q}=\hat{F}_f \ast \Phi$. Now, if $A_d=\Phi$ then, by means of theorem
\ref{eksdiaf} we get that the nonlinear system matches the desired linear one under any initial conditions ${\bf y}_0$.
If $\gcd (A_d,\Phi)=\tilde{\Phi}\ne \de_0$ then, by means of theorem
\ref{eksdiaf} we get that the nonlinear system matches the desired linear one under any initial conditions ${\bf y}_0$, which satisfy the relation $\tilde{\Phi}{\bf y}_0=0$. All the above are valid for any polynomial
$\Phi \in {\bf F}_1$ and hence the first branch of the feedback law of theorem \ref{mainresult} has been proved.
\par
\noindent 
 {\bf Case II.} 
 In this, rather extreme case, we want to find a common factor $\Phi$ of the terms of 
 (\ref{momatcheq}), by leaving the polynomial $S$ totally free. Let us suppose, as before, that ${\bf F}_2$ is non void and that $\Phi \in {\bf F}_2$. Step 5 of the MM-Algorithm implies that $RV\cap GV \ne \emptyset$.
This condition ensures that there is at least one rule ${\bf q}$ such that all the remainders are eliminated
and there is either a subset of the set $({\cal L}\cup \overline{\cal M})|_{\bf q} $ or a subset of the set $(\overline{\cal L}\cup {\cal M})|_{\bf q} $ consisting from polynomials $P$, with no zero companion coefficients
and with the property:
either $\gcd( P\in ({\cal L}\cup \overline{\cal M})|_{\bf q})=\Phi\ne\de_0$ or
$\gcd( P\in (\overline{\cal L}\cup {\cal M})|_{\bf q})=\Phi\ne\de_0$. Since  
$A\in {\cal L}$$\Rightarrow$$A=A_1\ast \Phi$ and
$B\in {\cal M}$$\Rightarrow$$B=B_1\ast \Phi$, for some polynomials $A_1,B_1$, so 
$A-B\ast S$$=(A_1-B_1 \ast S) \ast \Phi$, which implies that 
$\Phi$ is a factor of
$A-B \ast S$, too. From all the above we conclude that $\Phi$ is a common factor of all the terms of $F_f|_{\bf q}$, i.e. $F_f|_{\bf q}=\hat{F}_f \ast \Phi$, independently from the values of $S$. Therefore, if $A_d=\Phi$ then 
by means of \ref{eksdiaf} the nonlinear system and the linear desired one, will have the same dynamic behaviour 
for any initial conditions and any value of $S$, thus we set $S=Q$, $Q$ an arbitrary linear polynomial.
If $\gcd(A_d,\Phi)=\tilde{\Phi}\ne\de_0$, then 
by means of \ref{eksdiaf} the nonlinear system and the linear desired one will have the same dynamic behaviour 
for any initial conditions ${\bf y}_0$, which satisfy the equation
$\tilde{\Phi}{\bf y}_0=0$
and any value of $S$, thus we set again $S=Q$, $Q$ an arbitrary linear polynomial.
All the above are valid for any polynomial
$\Phi \in {\bf F}_2$ and hence the second branch of the feedback law of theorem \ref{mainresult} has been proved.
\par
\noindent 
 {\bf Case III.} We are dealing now with a case similar to the second one.
The essential difference is that the linear parts of $F$, which do not have parametrical coefficients are not involved to the calculation of the common factors.  
Let us suppose, as before, that ${\bf F}_3$ is non void and that $\Phi \in {\bf F}_3$. Step 6 of the MM-Algorithm implies that $RV\cap MV \ne \emptyset$.
This condition ensures that there is at least one rule ${\bf q}$ such that all the remainders are eliminated
and there is either a subset of the set $({\cal L}^*\cup \overline{\cal M}^*)|_{\bf q} $ or a subset of the set $(\overline{\cal L}^*\cup {\cal M}^*)|_{\bf q} $, consisting from polynomials $P$, with no zero companion coefficients
and with the property:
either $\gcd( P\in ({\cal L}^*\cup \overline{\cal M}^*)|_{\bf q})=\Phi\ne\de_0$ or
$\gcd( P \in (\overline{\cal L}^*\cup {\cal M}^*)|_{\bf q})=\Phi\ne\de_0$.
These relations mean that $\Phi$ is a common factor of all the terms of $F_f|_{\bf q}$, but the term
$A-B\ast S$. We have to choose now $S$ properly so that this latter term is a multiplier of $\Phi$, too.
This means that $A-B \ast S=R \ast \Phi $, for some polynomial
$R$ and thus $B \ast S + R \ast \Phi =A $. This is the well-known Diophantine
equation among linear polynomials, which has a solution if $\gcd(B,\Phi) | A$. If $\gcd(B,\Phi)=\de_0$ the above assumption is true. If $\gcd(B,\Phi)=\hat{\Phi}\ne \de_0$, then we conclude that $\gcd(A,B,\Phi)\ne\de_0$ but this implies that $\gcd( P\in ({\cal L}\cup \overline{\cal M})|_{\bf q})=\Phi\ne\de_0$ or
$\gcd( P\in (\overline{\cal L}\cup {\cal M})|_{\bf q})=\Phi\ne\de_0$ which leads us back to case II. 
Now, it is known that if $S_0$ is an initial solution of the Diophantine equation, $S$ must have the general type:
$S=S_0+Q\ast \Phi$, $Q$ an arbitrary linear polynomial. These values of $S$ ensure that $\Phi$ is a common factor of $F_f|_{\bf q}$ and working as in the previous cases we can prove the last branch of the feedback-law design.
 In all of the above, causality is
accomplished via a proper choice of the polynomial $Q$.
\begin{exa}
Let us suppose that we have the nonlinear system:
\[ F=y(t)-2y(t-1)+2y(t-1)y(t-2)+4y(t-2)^2+\frac{1}{2}y(t-2)y(t-3)-2y(t-2)y(t-4)+\]
\[+\frac{1}{4}y(t-4)^2+3u(t-2)+16u(t-1)^2-18u(t-1)u(t-2)-2u(t-1)^2u(t-2)+\]
\[+5u(t-2)^2+2u(t-1)u(t-2)^2-\frac{1}{2}u(t-2)^3+u(t-1)^2u(t-3)-u(t-1)u(t-2)u(t-3)+\]
\[+\frac{1}{4}u(t-2)^2u(t-3)+y(t-1)y(t-2)u(t-1)+\frac{1}{4}y(t-2)^2u(t-1)+\frac{1}{2}y(t-1)y(t-3)u(t-1)-\]
\[-\frac{1}{16}y(t-3)^2u(t-1)-4y(t-1)y(t-2)u(t-2)+y(t-2)^2u(t-2)+\]
\[+2y(t-1)y(t-3)u(t-2)-\frac{1}{4}y(t-3)^2u(t-2)\]
or using $\de\e$-polynomials $F[y(t),u(t)]=0$, where $F$ is the $\de\e$-polynomial $A$, appeared in the example \ref{exaLFL}. We want to find a feedback law $u(t)=Sy(t)$ and some classes of desired systems $A_d$, which are matched by the the closed-loop system. To illustrate our approach we shall follow the $MM$-Algorithm in details. We shall follow the $MM$-Algorithm. In order to construct the sets ${\bf F_1,F_2,F_3}$ we have first to calculate the sets of rules $LV,GV,MV,DR,RV$. Using the results of the example \ref{exaLFL} we see that $R_{\de,0}=0$ and thus $DR={\bf R}^2$, which means that any value of the parameters $(w_{\lambda,u,h,k},s_{\lambda,\omega,h,k})$ satisfies the condition
$R_{\de,0}=0$. We use now the $GCD$-Values Algorithm to discover common factors of the members of ${\cal L}_{\de}$. Let us further suppose that we are in  step 1 of this algorithm and that $n=1$, this suggests that we shall work with only one polynomial and all the others will be " eliminated " by putting their companion coefficients equal to zero. Since $ccoef(L_{1,1,0})=\frac{1}{4}$ and $ccoef(A)=1$, $A=\de_0-2\de_1$, this cannot be happened. Let us now suppose that $n=2$ and we have again to work only with the polynomials $L_{1,1,0}$ and $A$ by eliminationg the companion coefficients of all the others. A first glance at the coefficients indicates that the substitution $w_{1,1,0,3}=0$, $w_{1,1,0,2}=-4$, $w_{1,1,0,1}=0$, eliminates most of them, but finally the polynomial $L_{6,1,0}$ will remain to the " game " with $ccoef(L_{6,1,0})=\frac{1}{2}$ and thus we have to go the $n=3$ case. Here, we have a result. We obtain three polynomials with nonzero companion coefficients by means of the following rule of substitution:
$NZ=\{ w_{1,1,0,3}=0,w_{1,1,0,2}=-4,w_{1,1,0,1}=0,w_{6,1,0,1}=0,w_{6,1,0,0}=4\}$. These polynomials are ${\cal A}=\{A,L_{1,1,0},L_{6,1,0}\}$. Now, we want to find a subset of $NZ$ such that the polynomials ${\cal A}|_{NZ}=\{ \de_0-2\de_1,-4\de_2+\de_4+w_{1,1,0,0}\de_0,4\de_0+\de_2\}$ have a common factor, but since $\gcd(\de_0-2\de_1,4\de_0+\de_1)=\de_0$ this cannot be done and thus $CF=\emptyset$. Working similarly, we explore all the different cases for various values of $n$ and we finally conclude that always $CF=\emptyset$ and so $GV=\emptyset$, which implies that ${\bf F}_1=\emptyset $. Let us now investigate the case ${\bf F}_2$. In order to get ${\bf F}_2\ne \emptyset $ all the linear polynomials , including $A$ and $B$, must have common factors, but since $\gcd(A,B)=\de_0$ this is impossible and hence ${\bf F}_2=\emptyset$, too. It remains only the third case. This situation is similar with the previous one but without involving the no formal polynomials $A$ and $B$. This will provide us with some degrees of freedom. Indeed, the application of the $GCD$-Values Algorithm will supply us with the following common factor for the polynomials of the sets ${\cal L}^*$ and $\overline{\cal M}^*$: $\Phi=-2\de_0+\de_1$. Since $R_{\e,h}=0,R_{\de,h^*}=0$ we finally get that $\Phi \in {\bf F}_3 \ne \emptyset$. 
It remains to specify the desired systems $A_d$, which can be matched by the feedback-laws given in theorem \ref{mainresult}. If $A_d=\Phi$, then we seek for polynomials $S$ and $R$ such that  
\[ \de_0-2\de_1=R \ast(-2\de_0+\de_1)+ 3 \de_2 \ast S \] 
Solving this Diophantine equation we get a first solution $R_0=-\frac{1}{2}\de_0+\frac{3}{4}\de_1$ and $S_0=-\frac{1}{4}\de_0$. Therefore, the feedback-law upon request is $u(t)=Sy(t)$, with
$S=-2\de_0+\de_1+Q\ast (-\frac{1}{4}\de_0)$, $Q$ arbitrary linear polynomial. This law matches the closed-loop system with $A_d$, for any set of initial conditions. If $A_d$ has a common factor with $\Phi$, then $A_d$ must be a multiplier of $\Phi$, since $\Phi$ is a prime linear polynomial we conclude that $A_d=t \ast \Phi$, $T$ a linear polynomial. In this case, the feedback-law $u(t)=Sy(t)$, $S$ as before, creates a closed-loop system which matches the desired system for any set of initial conditions $(y_0,y_1)$, which satisfy the relation $-2y_1+y_0=0$, (theorem \ref{mainresult}). Finally, we choose $Q$ properly so that causality to be satisfied. For instance, by taking $Q=-8\de_0+2\de_2$, we get $S=\de_1-\frac{1}{2}\de_2$ or
\[ u(t)=y(t-1)-\frac{1}{2}y(t-2) \]
\end{exa}

\begin{exa}
Let us consider a very simplified model for velocity control
of an aircraft, \cite{kn:carcanias},
\[ \dot{x}_1=x_2-f(x_1), \quad \dot{x}_2=-x_2+u, \quad y=x_1 \]
where $x_1$ is the velocity, $x_2$ is the engine thrust and
$f(x_1)=x_1^2-2x_1$ is the aerodynamic drag. After some simple
manipulations we can take the following input-output relation in 
continuous form:
\[ \ddot{y}+(2y-1)\dot{y}+y^2-2y=u \]
To discretize the above system we choose to use the Euler forward
difference approximations, i.e. 
$\dot{y}=\frac{y(t+1)-y(t)}{h}$,
$\ddot{y}=\frac{1}{h^2}[y(t)-2y(t+1)+y(t+2)]$ and $h=1$.
A straight substitution and a proper shifting
of the time delays, will give:
\begin{equation}\label{aeroplano}
-\frac{1}{3}y(t)+y(t-1)-\frac{2}{3}y(t-1)y(t-2)+\frac{1}{3}y^2(t-2)=
-\frac{1}{3}u(t-2)
\end{equation}
or using $\de$ and $\e$-polynomials, $F[y(t),u(t)]=0$ with
$F=\frac{1}{3}\de_0+\de_1-\frac{2}{3}\de_1\de_2+\frac{1}{3}\de_2^2$,
$+\frac{1}{3}\e_2$. We want to design a feedback law
so that the closed-loop nonlinear system will match the stable linear
system: $6y(t)=y(t-1)+y(t-2)$ or
$(-6\de_0+\de_1+\de_2)y(t)=0$. We shall follow the MM-algorithm.
The absence of nonlinear $\e$-terms, as well as of $\de\e$-terms, simplifies the complexity of the operations. Step 2 will give the set:
${\cal L}=\{w_{10}\de_0+w_{11}\de_1+\de_2,$
$w_{20}\de_0+\de_1,$
$w_{30}\de_0+\de_1,$
$\frac{1}{3} \left( (-w_{10}^2+w_{11}^2w_{30}^2-2w_{20}w_{30}^2-2
w_{11}w_{20}w_{30}^2)\de_0+\right.$
$(-2w_{10}w_{11}+2w_{20}^2+2w_{11}w_{20}^2+$
$2w_{11}^2w_{30}-4w_{20}w_{30}-4w_{11}w_{20}w_{30})\de_1+$
$\left.(-2w_{10}+2w_{11}w_{20}+2w_{20})\de_2 \right)\}$.
At step 3 we form the following sets: $LV=\emptyset$,
$GV=MV=$
$\{w_{10}=-2x,w_{11}=x-2,$
$w_{20}=-2,w_{30}=-2, x \in {\bf R} \}$, (we did not include all
the operations in details for the sake of presentation) and we find a common factor
$\Phi=-2\de_0+\de_1$. Since $gcd(A_d,\Phi)=\Phi\ne\de_0$, theorem \ref{mainresult}
will give $S=S_0+Q \ast \Phi$, where $Q$ arbitrary, 
and $S_0$ satisfies the equation
$R \ast \Phi+S\ast B=A$ or
$R \ast (-2\de_0+\de_1) +$
$S \ast (-\frac{1}{3} \de_2)$
$=-\frac{1}{3}\de_0+\de_1$. The latter equation will give
$S=-\frac{5}{4}\de_0$ and
thus $S=-\frac{5}{4}\de_0+Q\ast (-2\de_0+\de_1)$.
In order to obtain a causal feedback connection we select
$Q=-\frac{5}{8}\de_0+\de_1$ and finally $S=-\frac{37}{8}\de_1+2\de_2$.
Hence, a feedback-law upon request is $u(t)=-\frac{37}{8}y(t-1)+2y(t-2)$.
This will law will produce a close-loop system wich matches $A_d$ for any set of initial conditions $(y_0,y_1)$ which satisfy the relation $-2y_1+y_0$. Please note that since the desired system is stable we have actually stabilized
the nonlinear system, by using a feedback which is stable as well.

\end{exa}

\section{Concluding Remark}

The aim of this paper was to describe an algebraic computational method for the study of
a general class of non-linear discrete input - output systems that contain products
between input and output signals. The model matching problem of these systems is
solved through certain symbolic algorithms. The entire approach is based on a
proper framework that involves the so-called $\de\e$-operators and the star-product
operation. We hope to be able to present current work on further applications of this
method to other questions in a future paper.

\end{document}